%% file: paper.tex
\documentclass[preprint]{elsarticle}

\usepackage{graphicx}
\usepackage{subfigure}
\usepackage{color}
\usepackage{amsmath}
\usepackage{algorithm}
\usepackage{amssymb}
\usepackage{latexsym}
\usepackage{amscd}
\usepackage{stmaryrd}
\usepackage{url}
\usepackage{natbib}
\usepackage{import}
\usepackage{booktabs}

\newdefinition{rmk}{Remark}

\newcommand{\M}{\mathcal{M}}

\newcommand{\dx}{\;\textrm{d}x}
\newcommand{\ds}{\;\textrm{d}s}
\newcommand{\dt}{\;\textrm{d}t}

\newcommand{\plotwidth}{10cm}
\newcommand{\figwidth}{12cm}

\newcommand{\R}{\mathcal{R}}

\newcommand{\uu}{\tilde{u}}
\newcommand{\vv}{\tilde{v}}
\newcommand{\UU}{\tilde{U}}
\newcommand{\zz}{\tilde{z}}
\newcommand{\ZZ}{\tilde{Z}}

\newcommand{\WW}{\hat{W}}
\newcommand{\ee}{\tilde{e}}

\newcommand{\apost}{\emph{a~posteriori}}

\newcommand{\RR}{\mathbb{R}}
\newcommand{\inner}[2]{\langle #1 ,\; #2 \rangle}
\newcommand{\renni}[2]{\langle #2 ,\; #1 \rangle}
\newcommand{\dinner}[2]{\langle\!\langle\!\langle #1 ,\; #2 \rangle\!\rangle\!\rangle}
\newcommand{\dnorm}[1]{|\!|\!|#1|\!|\!|}

\newcommand{\ipcsII}{
  \begin{algorithm}[htbp!]
    Let $k_n = t_n - t_{n-1}$ denote the time step and $I_n =
    (t_{n-1},t_n]$ the corresponding time interval. Furthermore, let
    $V_h$ and $\hat{V}_h$ denote a pair of trial and test spaces on a
    domain $\Omega$. For each time interval $I_n$, we seek the fluid
    velocity $U^n = U(\cdot, t_n)\in V_h$ and pressure $P^n = P(\cdot, t_n)\in Q_h$
    at time $t_n$ by solving the following three variational
    problems:
    \begin{enumerate}
    \item[1)] Compute the tentative velocity $U^{\bigstar}$ by
      solving
      \begin{equation}
        \begin{split}
          \renni{v}{\mathrm{d}_t^n(U^{\bigstar})}
          + \renni{\epsilon(v)}{\sigma(U^{n -\frac{1}{2}}, P^{n-1})} \\
          - \renni{v}{\nu (\nabla U^{n-\frac{1}{2}})^\top \, n}_{\Gamma_N}
          + \renni{v}{P^{n-1}n}_{\Gamma_N}
          = \renni{v}{f}
        \end{split}
      \end{equation}
      for all $v\in \hat{V}_h$, including any boundary conditions for
      the velocity. Here, $\mathrm{d}_t^n(U^{\bigstar}) =
      (U^{\bigstar} - U^{n-1}) / k_n + (U^{n-1}\cdot \nabla) U^{n-1}$
      and $U^{n-\frac{1}{2}} = (U^{\bigstar} + U^{n-1})/2$.\\
    \item[2)] Compute the corrected pressure $P^{n}$ by solving
      \begin{equation}
        \renni{\nabla q}{\nabla P^{n}}
        = \renni{\nabla q}{\nabla P^{n-1}}
        - k_n^{-1} \renni{q}{\nabla \cdot U^{\bigstar}}
      \end{equation}
      for all $q\in \hat{Q}_h$, including any boundary conditions for the pressure.\\
    \item[3)] Compute the corrected velocity $U^{n}$ by solving
      \begin{equation}
        \renni{v}{U^{n}}
        = \renni{v}{U^{\bigstar}}
        - k_n \renni{v}{\nabla(P^{n} - P^{n-1})}
      \end{equation}
      for all $v\in \hat{V}_h$,
      including any boundary conditions for the velocity.
    \end{enumerate}
    \caption{The Incremental Pressure Correction Scheme (IPCS)}
    \label{p2:alg:ipcs}
  \end{algorithm}
}

\begin{document}

\begin{frontmatter}

\title{An Adaptive Finite Element Splitting Method for the
Incompressible Navier--Stokes Equations}

\author[KS,focal]{K.~Selim\corref{cor1}\fnref{SIMULA} \fnref{UIO}}
\author[AL,focal]{A.~Logg\corref{cor2}\fnref{SIMULA} \fnref{UIO}}
\ead{logg@simula.no}
\author[ML]{M. G.~Larson\fnref{UMEA}}

\address[KS]{Simula Research Laboratory, P.O. Box 134, N-1325 Lysaker, Norway}
\address[AL]{Simula Research Laboratory, P.O. Box 134, N-1325 Lysaker, Norway}
\address[ML]{Department of Mathematics, Ume\aa\;University, SE-901 87 Ume\aa, Sweden}

\cortext[cor1]{Principle author}
\cortext[cor2]{Corresponding author}

\fntext[SIMULA]{Center for Biomedical Computing at Simula Research Laboratory, Norway}
\fntext[UIO]{ Department of Informatics, University of Oslo, Norway}
\fntext[UMEA]{Department of Mathematics, Ume\aa\;University, Sweden}

\begin{abstract}
  We present an adaptive finite element method for the incompressible
  Navier--Stokes equations based on a standard splitting scheme (the
  incremental pressure correction scheme). The presented method
  combines the efficiency and simplicity of a splitting method with
  the powerful framework offered by the finite element method for
  error analysis and adaptivity. An \apost{} error estimate is derived
  which expresses the error in a goal functional of interest as a sum
  of contributions from spatial discretization, time discretization
  and a term that measures the deviation of the splitting scheme from
  a pure Galerkin scheme (the computational error). Numerical examples
  are presented which demonstrate the performance of the adaptive
  algorithm and high quality efficiency indices. It is further
  demonstrated that the computational error of the Navier--Stokes
  momentum equation is linear in the size of the time step while the
  computational error of the continuity equation is quadratic in the
  size of the time step.
  \end{abstract}

\end{frontmatter}

\section{Introduction}

Adaptive finite element methods play an increasingly important role in
solving complex problems in science and engineering. The adaptive
methods are in general based on \apost{} error estimates, where the
error is estimated in terms of computable quantities, and adaptive
algorithms for automatic tuning of critical discretization parameters
such as the time step and the local mesh size.

Several important works have been published on \apost{} error
estimates for finite element approximations of time-dependent
problems, see for
example~\citet{ErikssonEstepEtAl1995,BeckerRannacher2001,
  GilesSuli2002}, the research monograph~\citet{EstepLarsonEtAl2000}
and the references therein. However, these estimates are generally
restricted to finite element approximations in space and time and thus
do not cover commonly used splitting schemes for efficient time
stepping. Splitting schemes are used to avoid solving coupled systems
of equations in each time step and have many applications, including
reaction--diffusion and fluid flow problems. An \apost{} error
estimate for a splitting method for systems of ordinary differential
equations was recently presented by~\citet{EstepGinting2008}. Error
analysis of non-Galerkin solutions is considered by
\citet{GilesSuli2002} with particular focus on error correction; that
is, improving the accuracy of a computed functional by
post-processing. \emph{A~posteriori} error analysis of the
incompressible Navier--Stokes has been studied in detail before, see
for example \citet{Hoffman2004,HoffmanJohnson2007}, but not for
splitting methods.

In this work, we consider splitting schemes for fluid flow. More
precisely, we derive an \apost{} error estimate for an incremental
pressure correction splitting scheme for the incompressible
Navier--Stokes equations. This type of scheme was originally proposed
by~\citet{Chorin1968} and \citet{Temam1969} and was later refined
by~\citet{Goda1979}. Even if a particular scheme is considered, our
approach extends to other splitting schemes or any other scheme.

The basic idea is to construct a piecewise polynomial interpolation in
time of the velocity and pressure and then apply the standard duality
argument to derive an error representation formula. Since the
splitting scheme does not satisfy a full Galerkin orthogonality, we
are left with an algebraic residual measuring the effect of the lack
of orthogonality caused by the splitting. The final estimate thus has
three contributions measuring the effect of discretization in space,
discretization in time and splitting, respectively. A similar approach
was briefly proposed but not implemented or tested
by~\citet{BengzonSplitting2009}. Based on the \apost{} error
estimates, we construct an adaptive algorithm and investigate the
performance of the adaptive algorithm and the quality of the error
estimate.

\subsection{Outline of this paper}

The outline of this paper is as follows. In the next section, we
present our model problem (the incompressible Navier--Stokes
equations). In Section~\ref{p2:sec:inconsistent-fem}, we introduce the
inconsistent finite element splitting method that is used to solve the
Navier--Stokes equations. The \apost{} error analysis is presented
in~Section~\ref{p2:sec:a-posteriori-analysis}, and the adaptive
algorithm is presented in Section~\ref{p2:sec:adaptive-algorithm}. The
efficiency of the adaptive algorithm and the quality of the error
estimate is demonstrated with numerical examples in
Section~\ref{p2:sec:numerical-results}. The paper closes with a
summary and some concluding remarks in
Section~\ref{p2:sec:conclusion}.

\section{The incompressible Navier--Stokes equations}
\label{p2:sec:preliminaries}

We consider a fluid governed by the incompressible Navier--Stokes
equations. For $\Omega \subset \RR^d$ $(d=2,3)$ we seek the velocity
$u:\Omega \times [0,T] \rightarrow \RR^d$ and pressure $p:\Omega
  \times [0,T] \rightarrow \RR$ such that
\begin{equation}
  \label{p2:eq:NS}
  \begin{array}{rlll}
     \dot{u} + (u\cdot \nabla)u - \nabla\cdot \sigma(u,p) & = & f
     &\quad \textrm{in}\;
    \Omega\times(0,T],\\
    \nabla\cdot u &=& 0 &\quad  \textrm{in}\;\Omega\times(0,T],
  \end{array}
\end{equation}
where $f$ is a given body force per unit volume.  In \eqref{p2:eq:NS},
the first equation is the momentum equation and the second equation is
the continuity equation. The symmetric Cauchy stress tensor
$\sigma(u,p)$ is defined as
\begin{equation}
  \sigma(u,p) = 2\nu\epsilon(u) - pI,
\end{equation}
where $\nu$ denotes the kinematic viscosity, $I$ is the identity
matrix and $\epsilon(u)$ is the symmetric gradient:
\begin{equation}
  \epsilon(u) = \tfrac{1}{2}(\nabla u + (\nabla u)^{\top}).
\end{equation}
We let $\Gamma = \Gamma_D \cup \Gamma_N$ denote the boundary of
$\Omega$ and associate Dirichlet and Neumann boundary conditions with
the two disjoint subsets $\Gamma_D$ and $\Gamma_N$, respectively. On
the Dirichlet boundary $\Gamma_D$, we impose a no-slip boundary
condition ($u = 0$) and assume a fully developed flow at the Neumann
boundary $\Gamma_N$; that is, $\nabla u \, n = 0$ where $n$ is the
outward pointing unit normal. This condition ensures that the flow
does not ``creep around the corners'' at the inflow and outflow. The
boundary condition is implemented weakly by dropping the term
involving $\nabla u$ from the boundary terms, leaving only $(\nu
(\nabla u)^{\top} - pI) \, n$.

\section{An inconsistent finite element formulation}
\label{p2:sec:inconsistent-fem}

Over the past couple of decades, numerous methods have been developed
for the numerical solution of the incompressible Navier--Stokes
equations. Many of these methods are based on a pseudo-compressibility
in order to overcome the algebraic difficulties of solving the
saddle-point problem resulting from a direct discretization of the
Navier--Stokes equations. A particular type of schemes are the
so-called splitting schemes where the velocity and pressure variables
are computed in a sequence of predictor--corrector type steps.

The first splitting method developed for the Navier--Stokes equations
is the Chorin projection scheme~\citep{Chorin1968,Temam1969}. Chorin's
scheme is a so-called non-incremental pressure correction scheme where
the starting point is to neglect the pressure in the momentum equation
and solve for a tentative velocity field. The tentative velocity is
then projected onto a divergence free space, resulting in a Poisson
problem for the pressure. In~\citet{Goda1979}, an improved version of
Chorin's scheme, the Incremental Pressure Correction Scheme (IPCS),
was presented, which solves the incompressible Navier--Stokes
equations in three steps. In the first step, an explicit pressure (the
value from the previous time step) is used in the momentum equation
and in the two subsequent steps, both the pressure and the velocity
are corrected.

In \citet{Valen-SendstadLoggEtAl2010a}, a comparison is made between
different splitting schemes, including a recent scheme
by~\citet{GuermondShen2003}, and (stabilized) Galerkin finite element
methods such as the G2 method by~\citet{HoffmanJohnson2007}. In this
study, six different numerical schemes were were tested on six
different test problems (making a total of 36 test cases). The test
cases all involved laminar flow at small to moderate size Reynolds
numbers in the range 1--1000. For each test problem, convergence in a
functional of interest or a global error norm was studied for a
sequence of refined meshes. The main conclusion of this study was that
the IPCS scheme was, overall, the most accurate and efficient method
for the particular choice of test problems. Based on these results, we
choose to base our implementation on the IPCS scheme, but emphasize
that the analysis is equally valid for any other scheme.

We consider the IPCS scheme in combination with a
Taylor--Hood~\citep{TaylorHood1973} approximation of the velocity and
pressure variables; that is, we seek a solution $\UU = (U,P)\in V_h \times
Q_h$, where $V_h$ is the space of continuous piecewise (vector-valued)
quadratic polynomials and $Q_h$ is the space of continuous piecewise
linear polynomials, respectively. A summary of the IPCS scheme is
given in Algorithm~\ref{p2:alg:ipcs}.

To analyze the error of the splitting method, one must construct a
suitable interpolant/continuous extension in order for the solution to
be defined at each point $(x, t) \in \Omega \times [0, T]$. Such an
interpolant comes natural for the scheme under consideration. Since
the solution is computed using a finite element formulation in space,
we only need to consider interpolation in time. In time, we define the
discrete solution to be the piecewise linear interpolant on each
interval $I_n$ based on the values $U^{n-1}, U^n$ and $P^{n-1}, P^n$,
respectively. For a higher order splitting scheme, care must be taken
in the construction of the interpolant to maintain the order of
accuracy.

\ipcsII

\section{\emph{A~posteriori} error analysis}
\label{p2:sec:a-posteriori-analysis}

To prove an \apost{} error estimate for the approximate solution
of~\eqref{p2:eq:NS} computed by the inconsistent finite element
formulation (the splitting scheme), we first state the weak form
of~\eqref{p2:eq:NS} in Section~\ref{p2:sec:weakprimal} and the
corresponding weak dual problem in Section~\ref{p2:sec:weakdual}.  We
then derive an error representation in
Section~\ref{p2:sec:errorrepresentation} from which we obtain the
error estimate(s) in Section~\ref{p2:sec:errorestimate}.

\subsection{The weak primal problem}
\label{p2:sec:weakprimal}

The weak form of~\eqref{p2:eq:NS} reads: find $(u,p)\in W = V\times Q =
\{v\in L^2(0,T;[H^1(\Omega)]^d):\: v(\cdot, 0) = u^0,
v|_{\Gamma_D}=0\} \times \{q\in L^2(\Omega):\: \langle q,\; 1\rangle
=0\}$ such that
\begin{equation} \label{p2:eq:weak-primal}
  a((u,p);(v,q)) = L((v,q))
\end{equation}
for all $(v,q)\in\WW$. The test space $\WW$ is defined analogously to
the trial space with homogeneous initial conditions.

The nonlinear form $a(\cdot;\cdot)$ and the linear form $L(\cdot)$
in~\eqref{p2:eq:weak-primal} are defined as
\begin{eqnarray}
  \nonumber
  a((u,p);(v,q)) &=& \int_0^T
  \renni{v}{\dot{u} + (u\cdot\nabla)u} + \renni{\epsilon(v)}{\sigma(u,p)} \\
  && \quad
  - \renni{v}{\nu(\nabla u)^\top \, n - pn}_{\Gamma_N}
  + \renni{q}{\nabla\cdot u} \dt, \label{p2:eq:adef}
  \\\nonumber\\
  L((v,q)) &=& \int_0^T \renni{v}{f}\dt.
\end{eqnarray}
We let $r : \WW \rightarrow \RR$ denote the weak residual
of~\eqref{p2:eq:weak-primal}; that is
\begin{equation} \label{p2:weak-residual}
  r((v, q)) = L((v,q)) - a((u,p);(v,q)) = \int_0^T r^t((v,q))\dt
\end{equation}
for all $(v,q)\in \WW$.

\subsection{The weak dual problem}
\label{p2:sec:weakdual}

Let now $\uu = (u,p)$ and let $\M = \M(\uu)$ denote a given linear
goal functional (the quantity of interest). The goal functional $\M$
is assumed to be of the form
\begin{equation} \label{p2:goal-functional}
  \M(\uu) = \M^T(u(\cdot, T)) + \int_0^T \M^t(\uu(\cdot,t))\dt.
\end{equation}
Here, $\M^T$ and $\M^t$ describe target functionals at the end
time $t = T$ and target functionals integrated over the time interval
$[0, T]$, respectively.

We may now introduce the weak dual problem of the incompressible
Navier--Stokes equations. The dual problem is used below in
Section~\ref{p2:sec:errorrepresentation} to express the error in the
goal functional $\M$ in terms of the weak
residual~\eqref{p2:weak-residual}. We let $\zz = (z,y)$ denote the
dual solution, where $z$ is the \emph{dual velocity} and $y$ the
\emph{dual pressure}. The (abstract) weak dual problem reads: find
$\zz\in W^*$ such that
\begin{equation} \label{p2:eq:weak-dual}
  \overline{a'}^*(\zz,\vv) = \M(\vv)
\end{equation}
for all $\vv \in \WW^*$. The dual trial and test spaces are defined by
$(W^*, \WW^*) = (\WW, W_0)$ where $W_0 = \{v-w : v,w \in W\}$. The
linearized, averaged and adjoint form $\overline{a'}^* : W^* \times \WW^*
\rightarrow \RR$ is defined by
\begin{equation} \label{p2:eq:linearization}
  \overline{a'}^*(\vv, \delta \uu)
  = \overline{a'}(\delta \uu, \vv)
  = \int_0^1 a'(s\uu + (1-s)\UU; \vv) \delta \uu \ds,
\end{equation}
where $a'$ denotes the Fr\'{e}chet derivative of the nonlinear form $a
: W \times \WW \rightarrow \mathbb{R}$ with respect to its first
argument.

To express the weak dual problem for the incompressible Navier--Stokes
equations, we start from the abstract dual
problem~\eqref{p2:eq:weak-dual} and differentiate the nonlinear form
$a$ defined in~\eqref{p2:eq:adef} with respect to the velocity field
$u$ and the pressure field $p$. The adjoint operator $^*$ amounts to
replacing the test functions $(v, q)$ in \eqref{p2:eq:adef} with the
dual variables $(z, y)$, and replacing the linearization variables
$\delta\uu = (\delta u, \delta p)$ in~\eqref{p2:eq:linearization} by
the dual test functions $(v, q)$. We find that the dual variational
problem reads: find $(z, y) \in W^*$ such that
\begin{equation} \label{p2:eq:fulldual}
  \begin{split}
    \int_0^T
    \inner{z}{\dot{v}} + \inner{z}{(\bar{u}\cdot\nabla)v +
    (v\cdot\nabla)\bar{u}} + \inner{\epsilon(z)}{\sigma(v,q)} \\
    - \inner{z}{\nu(\nabla v)^\top \, n - qn}_{\Gamma_N}
    + \inner{y}{\nabla\cdot v} \dt \\ = \M^T(v(\cdot, T))
    + \int_0^T \M^t((v, q))\dt
  \end{split}
\end{equation}
for all $(v, q) \in \hat{W}^*$, where $\bar{u} = \tfrac{1}{2}(U +
u)$. To solve~\eqref{p2:eq:fulldual}, we integrate the first term by
parts:
\begin{equation}
  \int_0^T \inner{z}{\dot{v}} \dt =
  \int_0^T \inner{-\dot{z}}{v} \dt +
  \inner{z(\cdot,T)}{v(\cdot,T)} - \inner{z(\cdot,0)}{v(\cdot,0)}.
\end{equation}
The boundary term at $t = 0$ vanishes since $(v, q) \in \WW^* = W_0$
and thus $v(\cdot, 0) = 0$. The second term cancels the term $\M^T(v(\cdot,
T))$ in the right-hand side of~\eqref{p2:eq:fulldual} if we take
$z(\cdot,T) = \psi^T$ where $\psi^T$ is the ($L^2$) Riesz
representer of $\M^T$. It follows that the dual solution may be
computed by solving the backward initial value problem
\begin{equation} \label{p2:eq:fulldual,2}
  \begin{split}
    \int_0^T
    \inner{-\dot{z}}{v} + \inner{z}{(\bar{u}\cdot\nabla)v +
    (v\cdot\nabla)\bar{u}} + \inner{\epsilon(z)}{\sigma(v,q)} \\
    - \inner{z}{\nu(\nabla v)^\top \, n - qn}_{\Gamma_N}
    + \inner{y}{\nabla\cdot v} \dt \\
    = \int_0^T \M^t((v, q))\dt,
  \end{split}
\end{equation}
with initial condition $z(\cdot,T) = \psi^T$.

\begin{rmk}
  The dual solution may be computed by a direct finite element
  discretization of \eqref{p2:eq:fulldual,2} with $\bar{u} \approx
  U$. In particular, it is not necessary to integrate by parts the
  remaining terms of \eqref{p2:eq:fulldual,2} to move derivatives from
  the test function~$v$.
\end{rmk}

\subsection{Error representation}
\label{p2:sec:errorrepresentation}

To derive a representation of the error $\M(\uu) - \M(\UU) = \M(\ee)$
in terms of the solution $\zz = (z, y)$ of the dual
problem~\eqref{p2:eq:weak-dual} and the weak residual~$r$ defined
in~\eqref{p2:weak-residual}, we note that by the definition of the
averaged linearized operator $\overline{a'}$
in~\eqref{p2:eq:linearization}, it follows that
\begin{equation}
  \begin{split}
    \overline{a'}(\ee, \vv)
    &= \int_0^1 a'(s\uu + (1-s)\UU; \vv) \ee \ds
    = \int_0^1 \frac{\mathrm{d}}{\mathrm{d}s} a(s\uu + (1-s)\UU; \vv) \ds \\
    &= a(\uu; \vv) - a(\UU; \vv),
  \end{split}
\end{equation}
for all $\vv \in W^*$, where $\ee = \uu- \UU \in \WW^*$. The error
representation now follows directly by taking $\vv = \ee \in \WW^*$
in~\eqref{p2:eq:weak-dual}:
\begin{equation} \label{p2:eq:error-representation}
  \M(\ee)
  = \overline{a'}^*(\zz, \ee)
  = \overline{a'}(\ee,\zz)
  = a(\uu; \zz) - a(\UU;\zz)
  = L(\zz) - a(\uu; \zz)
  = r(\zz).
\end{equation}
In other words, the error in the goal function $\M$ is the (weak)
residual of the dual solution.

If now the solution $\UU$ satisfies the Galerkin orthogonality $r(\vv)
= 0$ for all $\vv \in \WW_{hk} \subset \WW$, one may subtract a test
space interpolant $\pi_{hk}\zz$ to obtain $\M(\ee) = r(\zz
- \pi_{hk}\zz)$ from which the error estimate follows;
see~\citet{ErikssonJohnson1995,BeckerRannacher2001,BangerthRannacher2003}.
However, if the solution does not satisfy the Galerkin orthogonality,
one must account for the lack of orthogonality by adding and
subtracting the orthogonality condition. We do this in two steps to
separately account for the effects of space discretization, time
discretization and lack of orthogonality:
\begin{equation} \label{p2:eq:error-representation,2}
  \begin{split}
    \eta &\equiv \M(\ee) = r(\zz) \\
    &= r(\zz- \pi_h\zz + \pi_h \zz - \pi_{hk}\zz + \pi_{hk}\zz) \\
    &= r(\zz - \pi_h \zz) + r(\pi_h \zz - \pi_{hk}\zz) +  r(\pi_{hk}\zz) \\
    &\equiv \eta_h + \eta_k + \eta_c.
  \end{split}
\end{equation}
Here, $\pi_h$ is an interpolant into the semi-discrete space of
continuous piecewise quadratic vector-valued velocity fields and
continuous piecewise linear scalar pressure fields at each time $t\in
[0,T]$, $\pi_k$ is an interpolant into the semi-discrete space of
discontinuous piecewise constant functions (at each point $x \in
\Omega$), and $\pi_{hk} = \pi_k \pi_h$ is an interpolant into the
fully discrete test space $\WW_{hk} \subset \WW$.

\subsection{Error estimates}
\label{p2:sec:errorestimate}

To construct an adaptive algorithm based on the error
representation~\eqref{p2:eq:error-representation,2}, we estimate
$\eta_h$, $\eta_k$ and $\eta_c$ in terms of computable quantities to
obtain the total error estimate
\begin{equation}
  |\eta|
  = |\eta_h + \eta_k + \eta_c|
  \leq |\eta_h| + |\eta_k| + |\eta_c|
  \leq E_h + E_k + E_c
  \equiv E,
\end{equation}
where $|\eta_h| \leq E_h$, $|\eta_k| \leq E_k$ and $|\eta_c| \leq
E_c$.

\subsubsection{The space discretization error estimate $E_h$}
\label{p2:sec:space-error}

Starting from the definition $\eta_h = r(\zz -\pi_{h}\zz$), we
integrate by parts on each cell $K\in\mathcal{T}$, where $\mathcal{T}$
denotes the triangulation of $\Omega$, to obtain
\begin{equation} \label{p2:eq:error-estimate-h}
  \eta_h \leq \int_0^T \sum_{K\in\mathcal{T}} \eta_K \dt \equiv E_h,
\end{equation}
where $\eta_K = |\eta_K^1| + |\eta_K^2| + |\eta_K^3| + |\eta_K^4|$ and
\begin{eqnarray}
  \eta_K^1 &=&
  \renni{z - \pi_h z}{\dot{U} +  (U\cdot\nabla)U - \nabla \cdot \sigma(U, P) - f}_K, \\
  \eta_K^2 &=&
  \renni{z - \pi_h z}{\tfrac{1}{2}\llbracket\sigma(U, P)\rrbracket_n}_{\partial K \setminus \partial\Omega}, \\
  \eta_K^3 &=&
  \renni{z - \pi_h z}{\nu\nabla U \, n}_{\partial K \cap \Gamma_N}, \\
  \eta_K^4 &=&
  \renni{y - \pi_h y}{\nabla \cdot U}_K.
\end{eqnarray}
Here, $\llbracket\sigma(U,P)\rrbracket_n = \sigma(U^+,P^+) \, n^+ +
\sigma(U^-,P^-) \, n^-$ denotes the jump of the discrete normal stress
$\sigma(U,P) \, n$ across (interior) edges $\partial K$. The time
integral in~\eqref{p2:eq:error-estimate-h} is evaluated using midpoint
quadrature on each time interval $I_n$.

In practice, we approximate the dual solution $\zz$ by a numerical
approximation~$\ZZ$. However, care must be taken when inserting the
approximation~$\ZZ$ into the error
representation~\eqref{p2:eq:error-representation,2} or the error
estimate~\eqref{p2:eq:error-estimate-h}. In particular, the error
representation will evaluate to zero if the primal solution satisfies
the Galerkin orthogonality and the approximate dual solution is
computed on the same mesh and using the same order as the primal
solution. Furthermore, the error estimate will evaluate to zero since
$\ZZ - \pi_h \ZZ = \ZZ - \ZZ = 0$. Instead, we
compute an enhanced version $\mathcal{E}_h \ZZ$ from the computed
dual solution $\ZZ$ by local extrapolation on patches, as
described in~\citet{RognesLogg2010}. This allows the dual problem to
be solved on the same mesh using the same order as the primal problem,
which has many practical advantages.

\subsubsection{The time discretization error estimate $E_k$}
\label{p2:sec:time-error}

The time discretization error $\eta_k$ is estimated by
\begin{eqnarray}
\nonumber
  |\eta_k| &=&
  |r(\pi_h \zz - \pi_{hk} \zz)|
  = \left| \int_0^T r^t(\pi_h \zz - \pi_{hk} \zz) \dt \right|
  \leq \int_0^T |r^t(\pi_h \zz - \pi_{hk} \zz)| \dt\\
  &\equiv& E_k.
\end{eqnarray}
To evaluate the estimate $E_k$, we face the problem of integrating the
residual over each time interval $I_n$. This is challenging as the
residual oscillates heavily on each interval. For time discretizations
defined by a continuous Galerkin finite element method in time, the
residual is orthogonal to a space of discontinuous piecewise
polynomial functions on the partition of the time interval $[0,T]$. As
a consequence, the residual behaves like a Legendre polynomial on each
time interval~\citep{Logg2003a}. This is not necessarily the case for
a solution computed by a splitting method, as is the case
here. However, for the sake of analysis, we make the assumption that
the residual takes its maximum value at the endpoints of each interval
$I_n$. For a piecewise linear finite element approximation in time,
the corresponding test space consists of the space of discontinuous
piecewise constant functions. We may then take the interpolant $\pi_k$
to be the midpoint value on each interval to obtain the estimate
\begin{equation} \label{p2:eq:ek1}
  \begin{split}
    E_k
    &\leq \sum_{n=1}^M k_n |r^t(\ZZ(\cdot, t_n)) -
    r^t((\ZZ(\cdot, t_{n-1}) + \ZZ(\cdot, t_n))/2)| \\
    &= \frac{1}{2}
    \sum_{n=1}^M k_n |r^t(\ZZ(\cdot, t_n)) - r^t(\ZZ(\cdot, t_{n-1}))|,
  \end{split}
\end{equation}
where $\ZZ$ is the approximate numerical solution of the dual
problem and $M$ is the number of time steps.

The estimate~\eqref{p2:eq:ek1} is used to estimate the size of the
time discretization error $\eta_k$ by the adaptive algorithm presented
below in Section~\ref{p2:sec:adaptive-algorithm}. To control the size
of the adaptive time steps, we here derive an alternate estimate. We
let $R^t$ denote the $L^2$ Riesz representer of the functional $r^t$
and write
\begin{equation} \label{eq:Ek,S}
  \begin{split}
    E_k
    &= \int_0^T |r^t(\pi_h \zz - \pi_{hk} \zz)| \dt
    = \int_0^T |\dinner{\pi_h \zz - \pi_{hk} \zz}{R^t} | \dt \\
    &\leq \int_0^T \dnorm{\pi_h \zz - \pi_{hk} \zz} \; \dnorm{R^t} \dt \\
    &\leq
    \max_{[0,T]}\{k_n(t) \dnorm{R^t} \}
    \int_0^T k_n^{-1} \dnorm{\pi_h \zz - \pi_{hk} \zz} \dt \\
    &= S(T) \max_{[0,T]}\{k_n(t) \dnorm{R^t} \} \\
    &\equiv \bar{E}_k,
  \end{split}
\end{equation}
where $S(T) = \int_0^T k_n^{-1} \dnorm{\pi_h \zz - \pi_{hk} \zz} \dt$
is a stability factor. The inner product $\dinner{\cdot}{\cdot}$ is
here defined by $\dinner{\uu}{\vv} = \inner{u}{v} + \inner{p}{q}$ and
the norm $\dnorm{\cdot}$ is defined as the corresponding norm. We
remark that the introduction of inequalities in~\eqref{eq:Ek,S} may
render the estimate less sharp. This however is not a problem
since~\eqref{eq:Ek,S} is not used to estimate the error; it is only
used to drive the selection of adaptive time steps.

The norm $\dnorm{R^t}$ of the Riesz representer may be
computed directly as follows. We first note that the Riesz representer
$R^t$ is defined by the variational problem
\begin{equation} \label{p2:riesz}
  \dinner{R^t}{\vv} = r^t(\vv)
\end{equation}
for all test functions $\vv \in \hat{W}_h$. The variational problem
\eqref{p2:riesz} corresponds to a linear system
\begin{equation} \label{p2:eq:riesz-system}
  M \R = b
\end{equation}
where $M$ is the mass matrix and $\R$ is the vector of degrees of
freedom for $R^t$. Clearly, the solution to~\eqref{p2:eq:riesz-system}
is given by $\R = M^{-1} b$. It follows that
\begin{equation}
  \begin{split}
    \dnorm{R^t}^2
    &= \dinner{R^t}{R^t}
     = \dinner{\sum_{i=1}^N \R_i \varphi_i}{\sum_{j=1}^N \R_j \varphi_j}
     = \sum_{i,j=1}^N \R_i \dinner{\varphi_i}{\varphi_j} \R_j \\
    &= (M^{-1}b)^\top M \R = b^\top \R,
  \end{split}
\end{equation}
by the symmetry of $M$. The residual norm $\dnorm{R^t}$ may thus be
computed by assembling and solving the linear
system~\eqref{p2:eq:riesz-system}, computing the inner product
$b^{\top}\R$ and taking the square root.

\subsubsection{The computational error estimate $E_c$}
\label{p2:sec:computational-error}

The computational error $\eta_c$ is computed by a direct evaluation
of the weak residual for the computed approximate dual solution
$\ZZ$:
\begin{equation} \label{p2:eq:Ec}
  |\eta_c| = |r(\pi_{hk}\zz)| \approx |r(\pi_k\ZZ)|
  = |\int_0^T r^t(\pi_k\ZZ) \dt | \equiv E_c,
\end{equation}
where $\ZZ$ is the approximate solution of the dual problem computed
on the same mesh and using the same polynomial degree and time steps
as the primal solution. In our implementation, we have chosen to
compute the dual solution by a simple application of the
$\mathrm{dG}(0)$ (backward Euler) method to the linear dual
problem. We note that, by construction, the computational error
estimate $E_c$ is zero if the primal solution satisfies the Galerkin
orthogonality.

\section{Adaptive algorithm}
\label{p2:sec:adaptive-algorithm}

Based on the \apost{} error estimate derived in
Section~\ref{p2:sec:errorestimate}, we may now formulate an adaptive
algorithm for the incompressible Navier--Stokes equations. The
adaptive algorithm is summarized in
Algorithm~\ref{p2:alg:adaptive-algorithm}.

\begin{algorithm}
  Given a goal functional $\M = \M(\uu)$ and a tolerance
  $\textrm{TOL}>0$:
  \begin{enumerate}
  \item[0)] Select an initial coarse mesh and initial time step.
  \item[1)] Solve the primal problem~\eqref{p2:eq:NS} using (for
    example) the incremental pressure correction scheme
    (Algorithm~\ref{p2:alg:ipcs}) on the current (fixed) mesh using
    adaptive time steps.
  \item[2)] Solve the dual problem~\eqref{p2:eq:fulldual,2} backward
    in time on the same mesh as the primal problem and using the same
    adaptive time steps.
  \item[3)] Evaluate the error estimate $E = E_h + E_k + E_c$ defined
    in~\eqref{p2:eq:error-estimate-h}, \eqref{p2:eq:ek1} and
    \eqref{p2:eq:Ec}, and the error indicators $\{\eta_K\}$.
  \item[4)] If $E \leq \mathrm{TOL}$, then stop.
  \item[5)] Refine the mesh based on the error indicators
    $\{\eta_K\}$.
  \item[6)] Continue from step 1).
  \end{enumerate}
  \caption{Adaptive algorithm}
  \label{p2:alg:adaptive-algorithm}
\end{algorithm}

In Algorithm~\ref{p2:alg:adaptive-algorithm}, we make use of two
different tolerances to $\textrm{TOL}_h$ and $\textrm{TOL}_k$ which
are used to control the errors in the space and time discretization
such that $\textrm{TOL}_h + \textrm{TOL}_k \leq \mathrm{TOL} -
E_c$. The computational error $E_c$ is only used as part of the
stopping criterion $E \leq \mathrm{TOL}$; it is not used to drive the
adaptive refinement. However, as will be demonstrated in
Section~\ref{p2:sec:numerical-results}, the computational error is
reduced when the size of the time step is reduced. One may therefore
consider extending the adaptive algorithm to control also the size of
the computational error $E_c$.

In each adaptive iteration, consisting of a full solution of the
primal problem, the dual problem and evaluation of the error
indicators, the mesh is adaptively refined based on fixed fraction
marking; that is, a fixed top fraction of the cells with the largest
indicators are marked for refinement. For mesh refinement, we consider
two different refinement strategies: the Rivara recursive bisection
algorithm~\citep{Rivara1991} and a regular cut algorithm which
subdivides all marked triangles into four congruent subtriangles and
propagates the refinement to neighboring triangles using bisection.

The step size $k_n$ is determined in each time step based on the error
estimate $\bar{E}_k = S(T) \max_{[0,T]}\{k_n(t) \dnorm{R^t}\}$. To
achieve $\bar{E}_k = \mathrm{TOL}_k$, we set
\begin{equation} \label{p2:eq:kn}
  k_n = \frac{\mathrm{TOL}_k}{S(T) \max_{[t_{n-1},t_n]} \dnorm{R^t}}
  = \frac{\mathrm{TOL}_k}{S(T) \dnorm{R^n}},
\end{equation}
where again we have made the assumption that the residual takes its
maximum value at the endpoints. Since $R^n$ is not known until the
solution has been computed on the time interval $I_n$, which in turn
depends on the size of the time step $k_n$, it is tempting to replace
$R^n$ by $R^{n-1}$ in~\eqref{p2:eq:kn}. However, this leads to
oscillations in the time step; if $R^{n-1}$ is large, $k_n$ will be
small and, as a consequence, $R^{n}$ will be small, which in turn
leads to a large step $k_n$ and so on. To control the time step,
one may introduce a form of smoothing by letting $\tilde{k}_n$ be the
time step determined by
\begin{equation}
  \tilde{k}_n = \frac{\mathrm{tol}_k}{\dnorm{R^{n-1}}},
\end{equation}
for $\mathrm{tol}_k = \mathrm{TOL}_k / S(T)$ and then take $k_n$ to be
the harmonic mean
\begin{equation}
  k_n = \frac{2 k_{n-1} \tilde{k}_n}{k_{n-1} + \tilde{k}_n}.
\end{equation}
See~\citet{Soderlind2002} and \citet{Logg2003b} for a further
discussion on time step selection. In practice, we do not compute the
stability factor $S(T)$ but instead adjust the size of
$\mathrm{tol}_k$ based on the size of $E_k$.

\section{Numerical results}
\label{p2:sec:numerical-results}

We here present numerical results to test the adaptive algorithm and
the quality of the derived error estimates. An implementation of the
adaptive solver, including the test problems described in this
section, is freely available as part of the open source solver
package~\citet{www:cbc.solve}. The package relies on the FEniCS/DOLFIN
finite element library~\citep{www:fenics,Logg2007a,LoggWells2010}.

\subsection{Case 1: Channel flow with wall-mounted body}

\begin{figure}[tbp!]
  \begin{center}
    \def\svgwidth{15cm}
    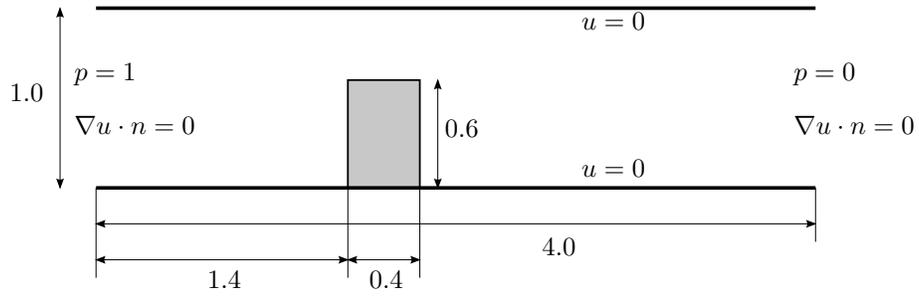
    \caption{Geometry and boundary conditions for the ``channel with
      flap'' model problem.}
    \label{p2:fig:geometry}
  \end{center}
\end{figure}

\begin{figure}[tbp!]
  \begin{center}
    \includegraphics[width=\figwidth]{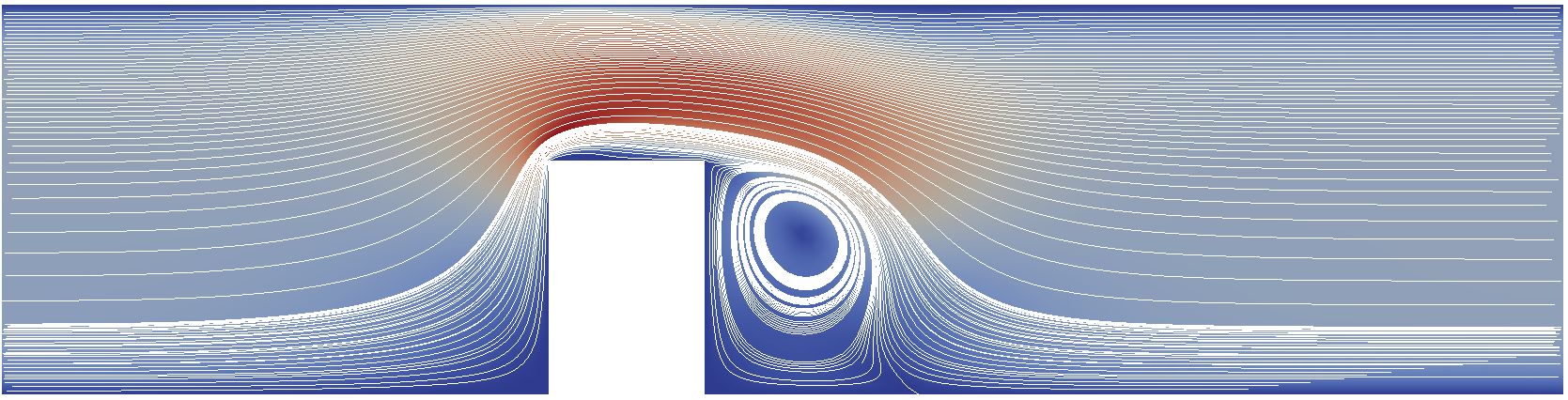}
    \includegraphics[width=\figwidth]{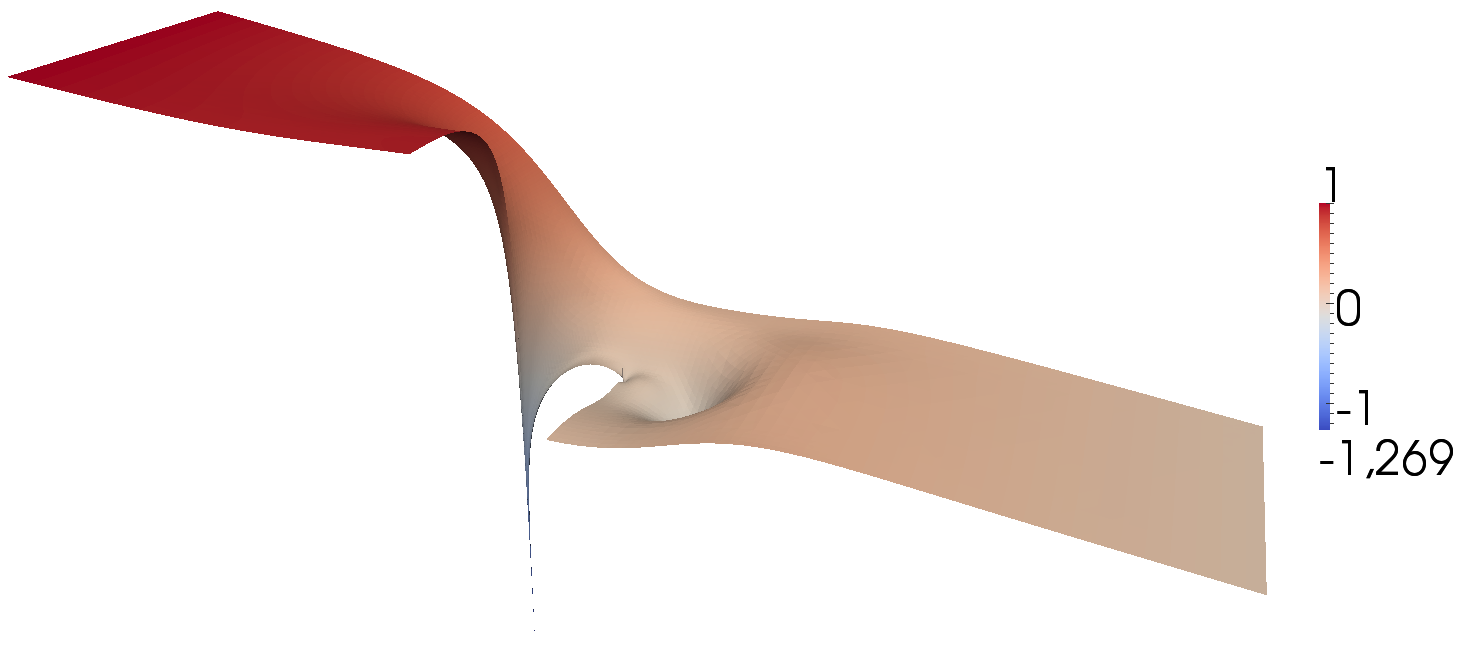}
    \caption{Fluid velocity (top) and pressure (bottom) at final time
      $T = 2.5$ computed with fixed time step $k = 0.005$ and
      14~levels of bisection refinement (marking fraction $0.3$). The
      final mesh has $16,581$ triangles ($76,085$ degrees of
      freedom). The colorbar indicates the range of the scalar
      pressure field.}
    \label{p2:fig:solution}
  \end{center}
\end{figure}

As a first test problem, we consider a wall-mounted body (a ``flap'')
immersed in a pressure-driven channel flow as illustrated in
Figure~\ref{p2:fig:geometry}. The kinematic viscosity is $\nu =
0.002$.  As initial condition, we set $u = 0$. The pressure boundary
conditions $p = 1$ at the inflow and $p = 0$ at the outflow accelerate
the flow from the initial stationary (zero) solution to the flow field
depicted in Figure~\ref{p2:fig:solution} at final time $T = 2.5$. Note
that the solution at final time is not stationary which is important
when measuring the performance and propagation of time discretization
errors.

As a goal functional, we consider the integrated shear stress on the
top of the flap:

\begin{equation} \label{p2:eq:goal0}
  \begin{split}
    \M_1(\uu)
    &= \int_0^T \int_{\Gamma_1} (\sigma(u, p) \, n) \cdot t \ds \dt \\
    &= \int_0^T \int_{\Gamma_1} \sigma_{12}(u, p) \ds \dt
    = \int_0^T \int_{\Gamma_1}
    \nu (\partial u_1 / \partial x_2 + \partial u_2 / \partial x_1) \ds \dt,
    \end{split}
\end{equation}
where $n = (0, 1)$, $t = (1, 0)$ and $\Gamma_1 = \{(x_1, x_2) : 1.4
\leq x_1 \leq 1.8, x_2 = 0.6\}$. As a reference value for the goal
functional, we take $\M_1(\uu) = 0.0200$. This reference value was
obtained by extrapolation from solutions computed with constant time
step $k = 0.005$ on a sequence of adaptively refined meshes.

\subsubsection{Dual solutions and adaptive meshes}

\begin{figure}[tbp!]
  \begin{center}
    \includegraphics[width=\plotwidth]{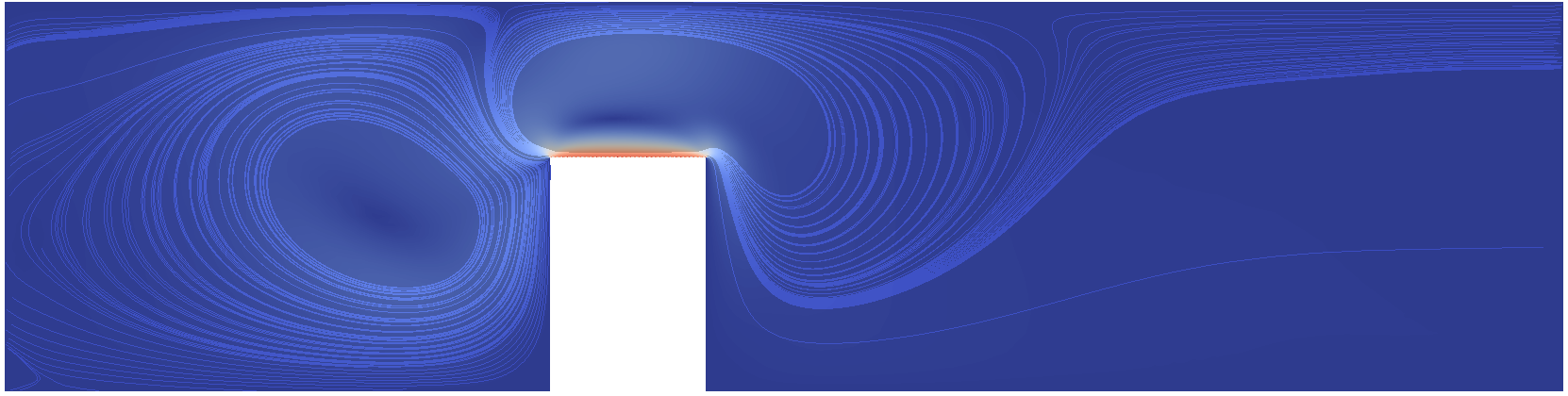}
    \includegraphics[width=\plotwidth]{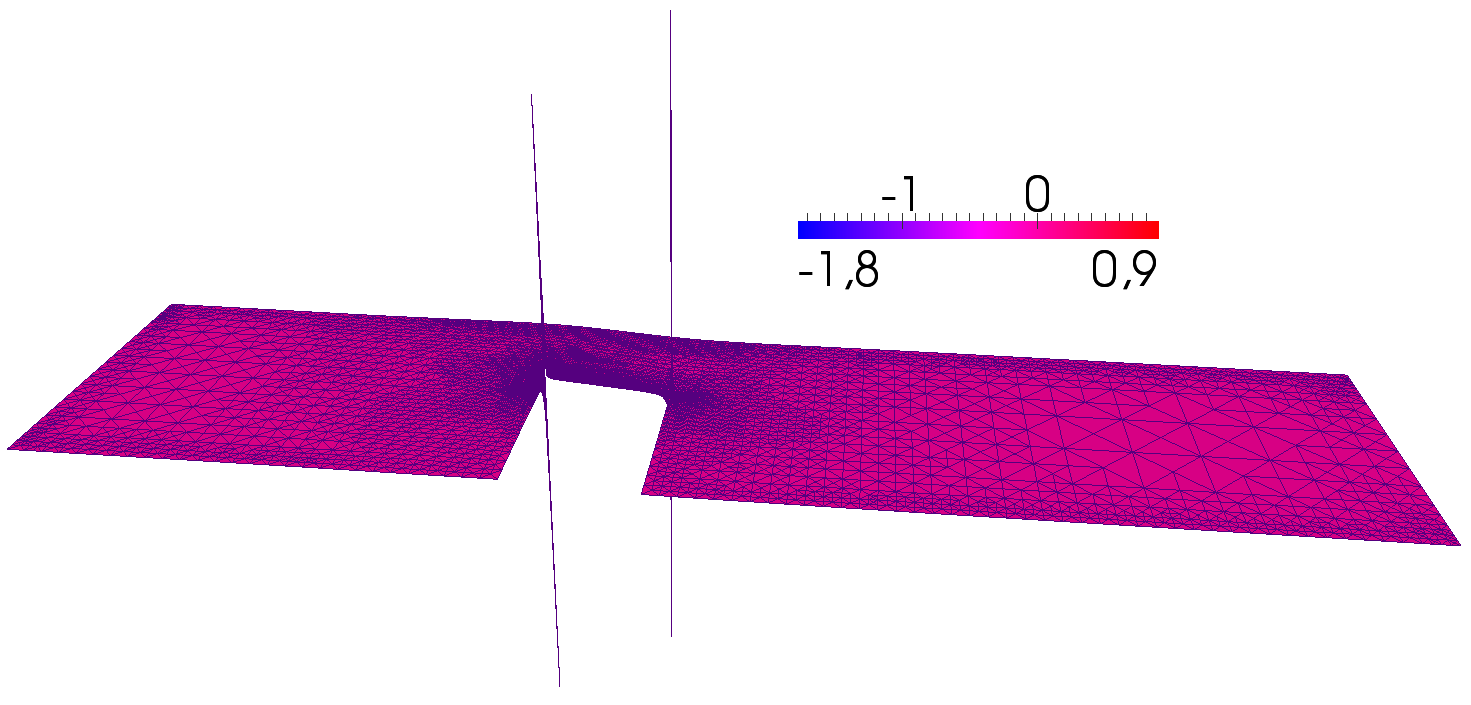}
    \caption{Dual fluid velocity (top) and dual pressure (bottom) at
      ``final'' time $t = 0$ for the channel flow test problem. The
      colorbar indicates the range of the scalar dual pressure field.}
    \label{p2:fig:dual0}
  \end{center}
\end{figure}

\begin{figure}[tbp!]
  \begin{center}
    \includegraphics[width=\plotwidth]{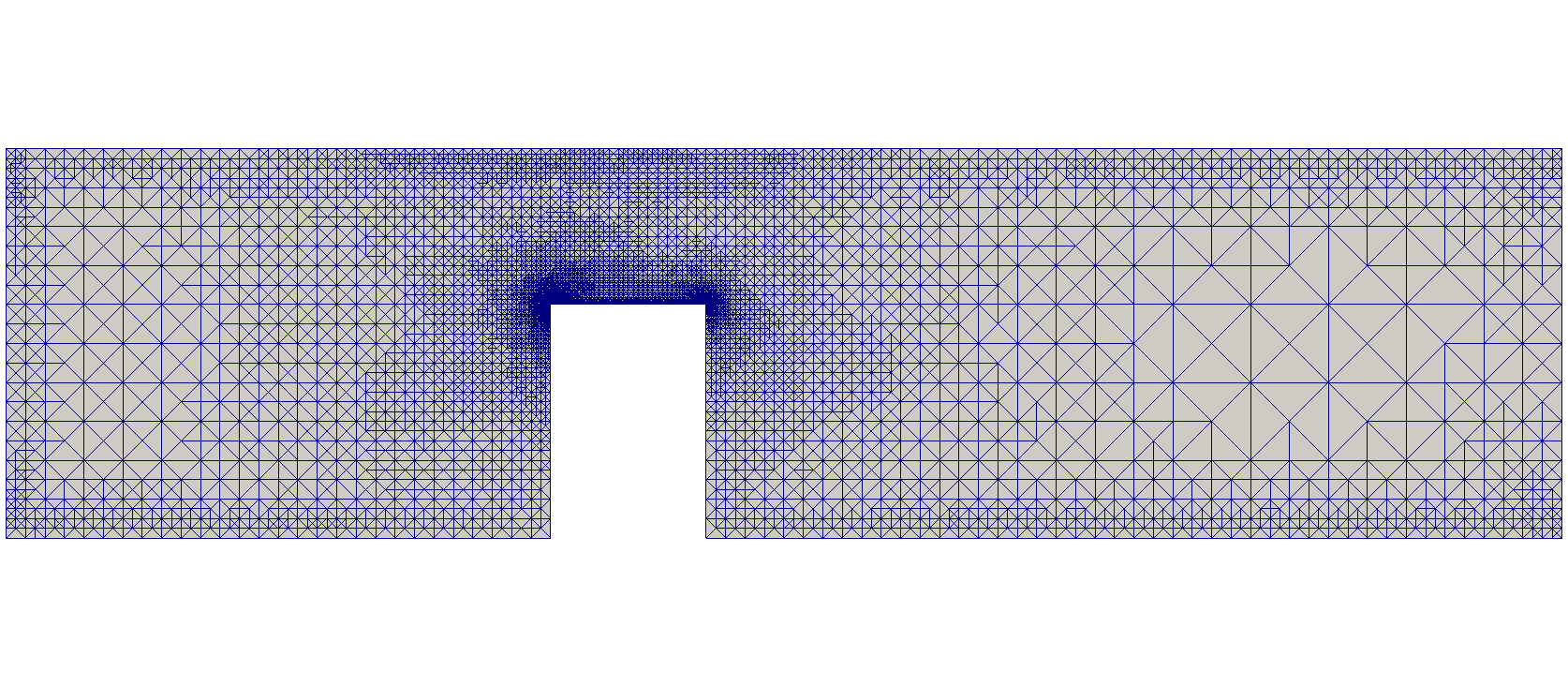}
    \caption{Mesh obtained by 14 levels of recursive bisection
      refinement with marking fraction $0.3$ for the channel flow test
      problem.}
    \label{p2:fig:mesh0}
  \end{center}
\end{figure}

The dual solutions corresponding to the goal functional $\M_1$ are
shown in Figure~\ref{p2:fig:dual0}. As seen in this Figure, the dual
solution clearly reflects the choice of goal functional. The dual
velocity~$z$ (and dual velocity gradients) are large close to the top
of the flap where the goal functional $\M_1$ measures the shear
stress. The same figure displays large spikes in the dual pressure~$y$
at the reentrant corners. A detail of the dual pressure spikes is
displayed in Figure~\ref{p2:fig:dual0,pressurezoom}. Combined with
large residuals in the vicinity of the reentrant corners, the result
is heavy refinement in a region located close to the top of the flap
as seen in Figure~\ref{p2:fig:mesh0}.

\begin{figure}[tbp!]
  \begin{center}
    \includegraphics[width=\plotwidth]{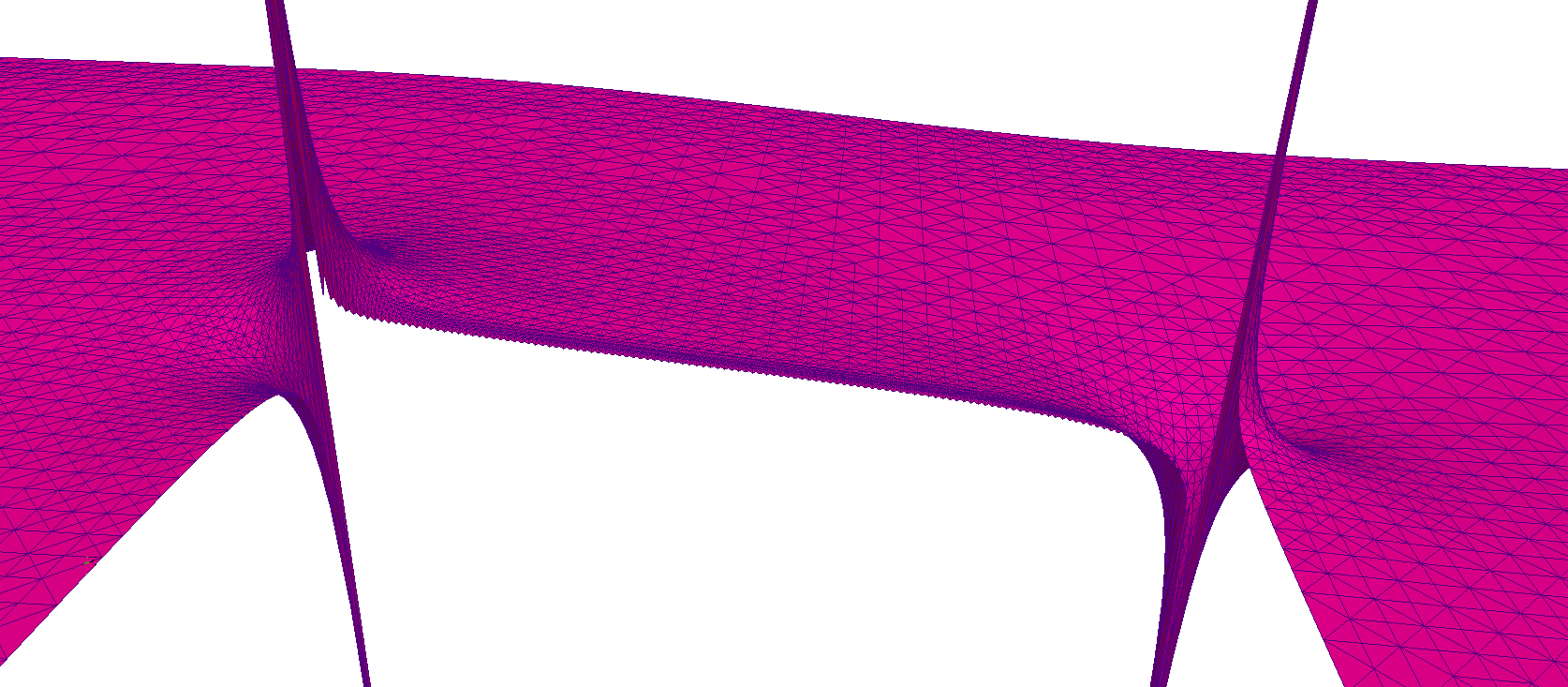}
    \caption{Detail of the dual pressure field at ``final'' time $t =
      0$ for the channel flow test problem.}
    \label{p2:fig:dual0,pressurezoom}
  \end{center}
\end{figure}

\subsubsection{Convergence and efficiency indices}

\begin{figure}[tbp!]
  \begin{center}
    \includegraphics[width=\plotwidth]{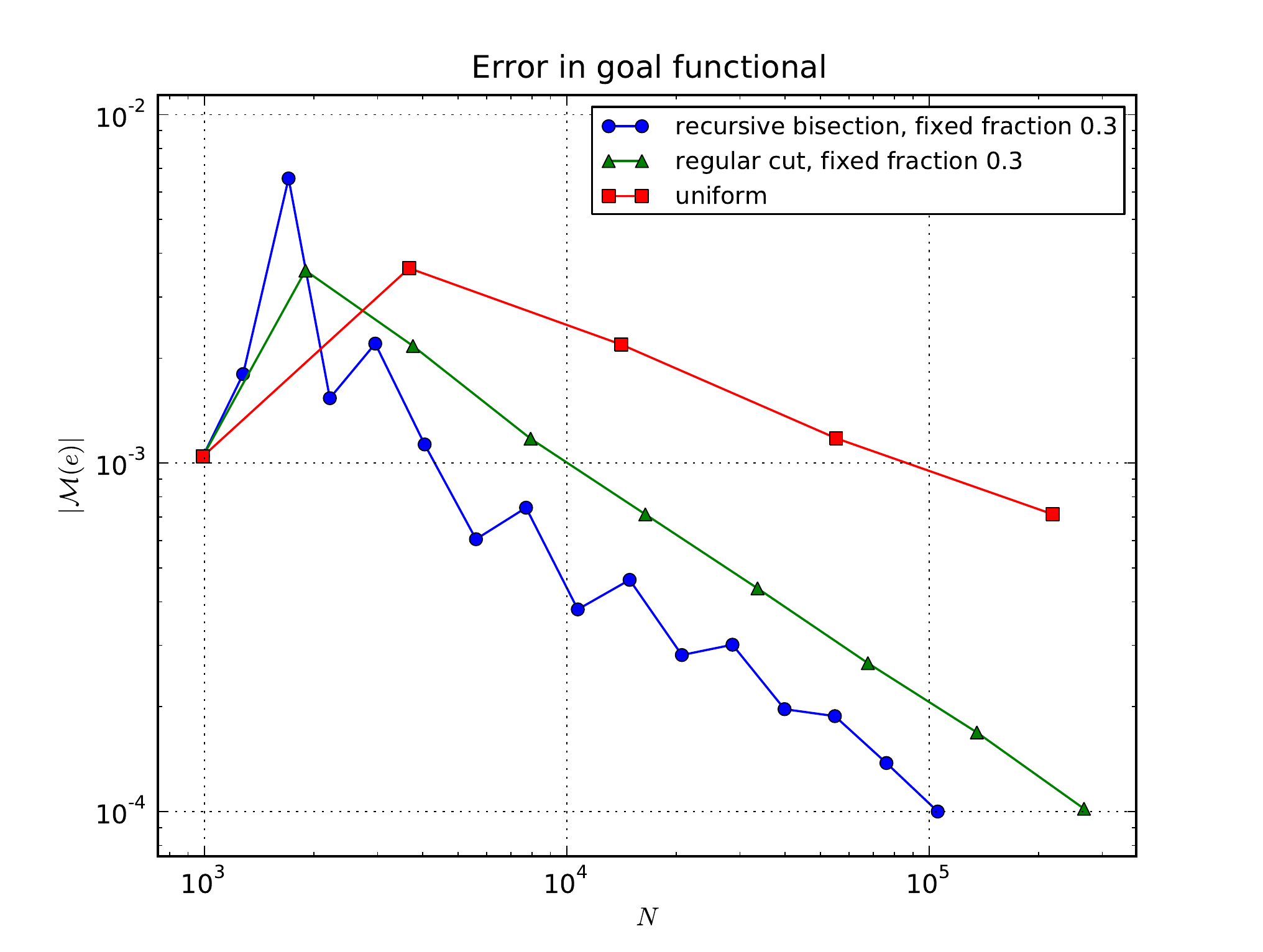}
    \includegraphics[width=\plotwidth]{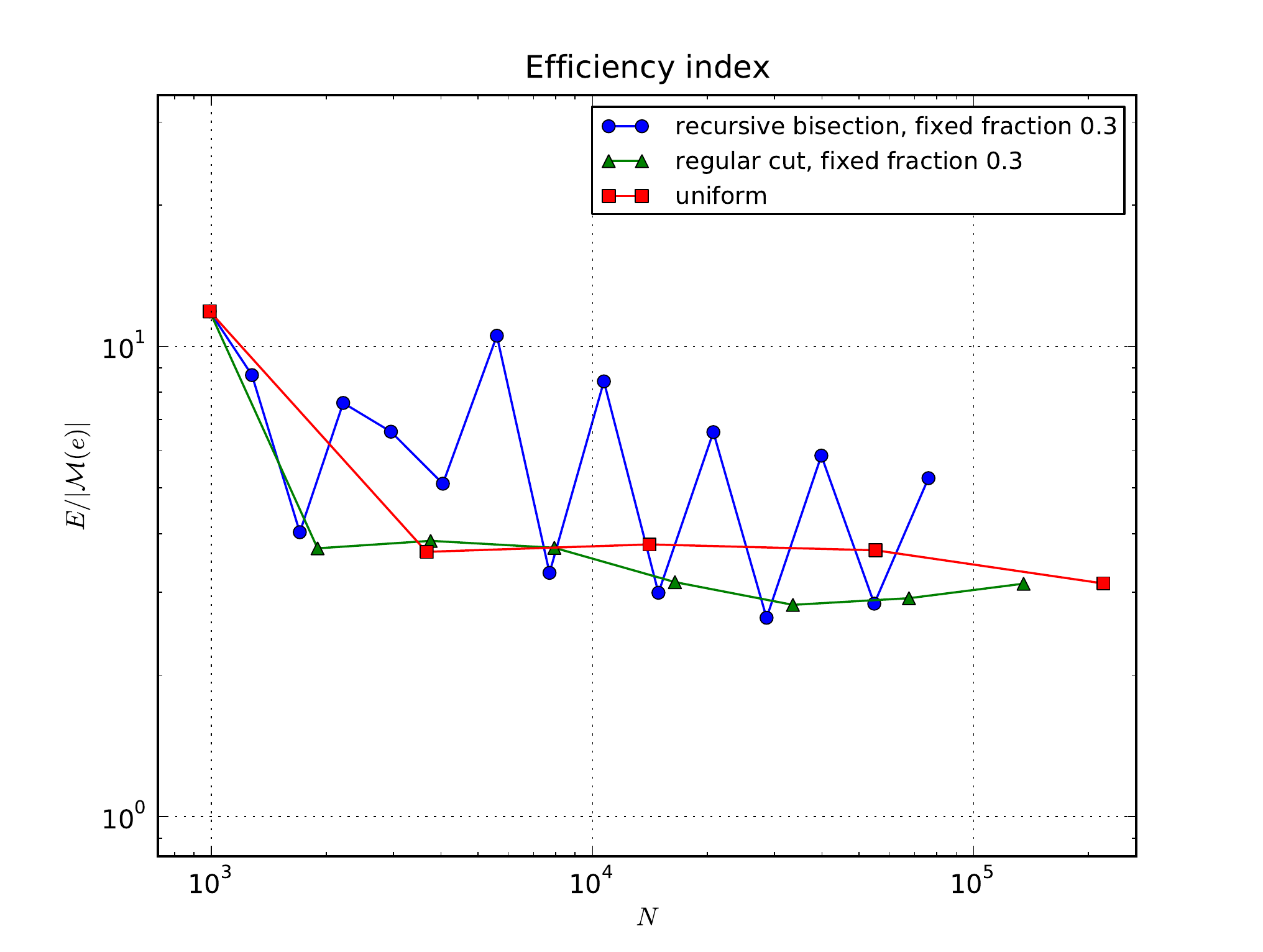}
    \caption{Error (top) and efficiency indices (bottom) as function
      of the number of spatial degrees of freedom for fixed time step
      $k = 0.005$, fixed fraction marking (marking fraction $0.3$) and
      varying refinement algorithms (recursive bisection, regular cut
      and uniform) for the channel flow test problem.}
    \label{p2:fig:conv0}
  \end{center}
\end{figure}

To investigate the performance of adaptive mesh refinement and the
quality of computed error estimates, we plot in
Figure~\ref{p2:fig:conv0} errors and efficiency indices for a sequence
of adaptively refined meshes and fixed time step $k = 0.005$. A
comparison is made between three different refinement algorithms:
recursive bisection, regular cut refinement and uniform (non-adaptive)
refinement. For both recursive bisection and regular cut refinement,
we use a fixed fraction marking strategy with marking fraction~$0.3$;
that is, in each refinement step, we mark for refinement the top
$30\%$ of all triangles with the largest error indicators.

We find that the adaptive algorithm performs very well; a uniformly
refined mesh requires more than an order of magnitude more degrees of
freedom to reach a prescribed tolerance. This is evident in
Figure~\ref{p2:fig:conv0} by finding the point where the error reaches
the level $|\M_1(\ee)| \leq 0.001$. This level is reached for roughly
$90,000$ degrees of freedom on a uniformly refined mesh, whereas the
adaptively refined meshes obtained by recursive bisection and regular
cut refinement reach the same level of accuracy using only $5,000$ and
$10,000$ degrees of freedom, respectively. We also note that while the
solution obtained by recursive bisection converges fastest, the
convergence of the solution obtained by regular cut refinement is more
robust. Computed efficiency indices (error estimate divided by actual
error) are stable and vary between ca. 3 and 4, which means that we
overestimate the error, but not by much.

\begin{figure}[tbp!]
  \begin{center}
    \includegraphics[width=\plotwidth]{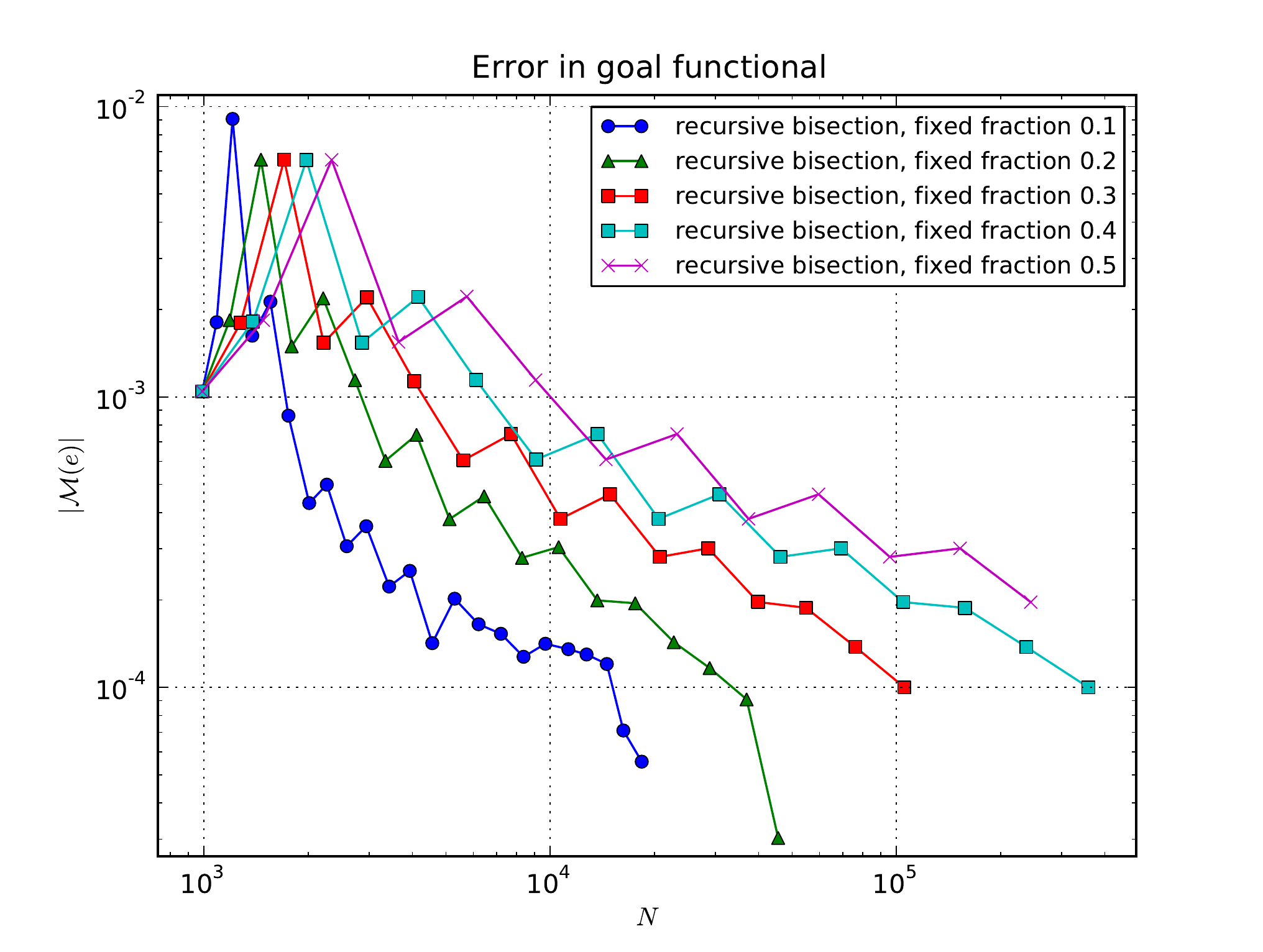}
    \includegraphics[width=\plotwidth]{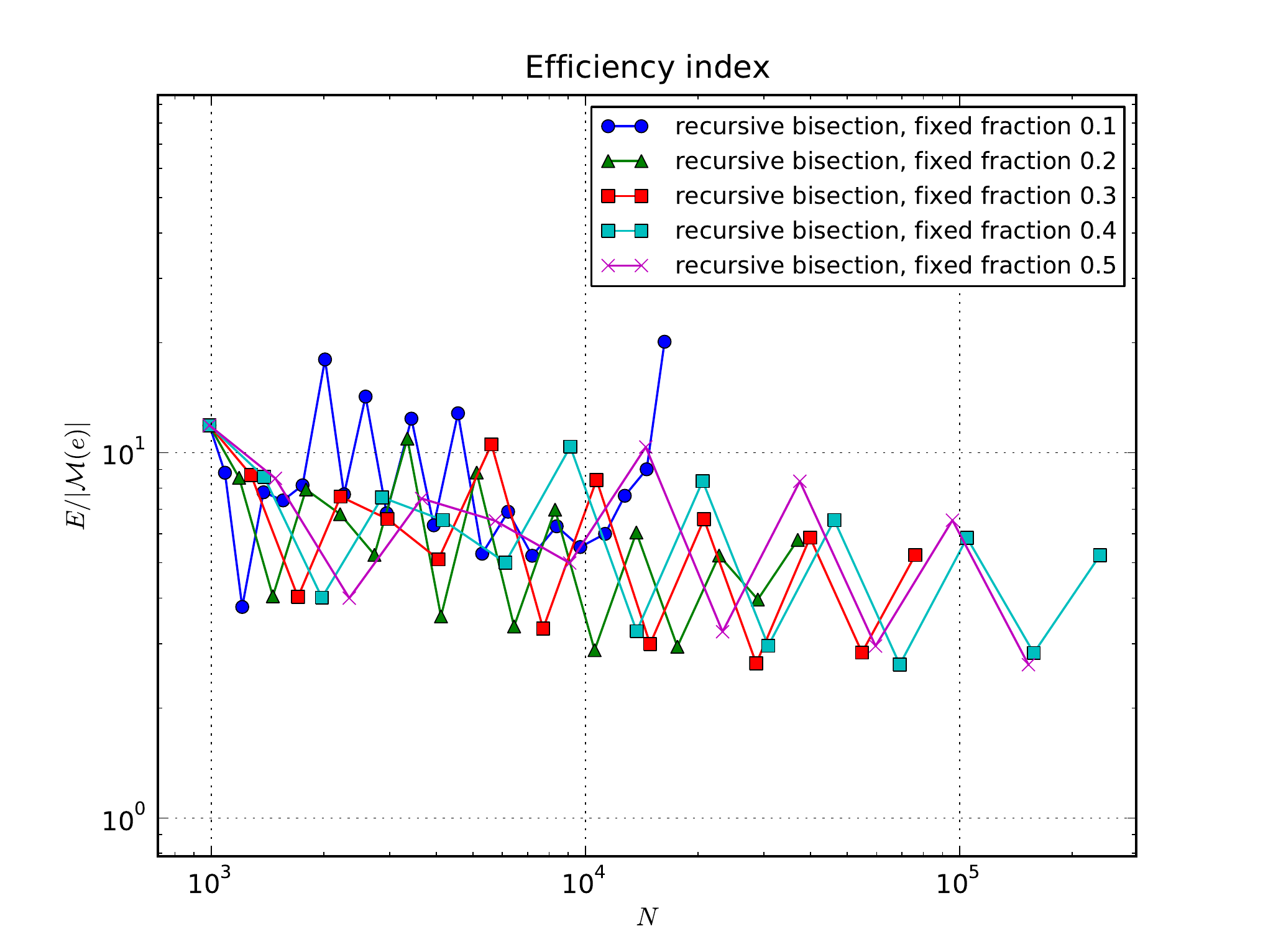}
    \caption{Comparison of errors (top) and efficiency indices
      (bottom) for varying marking fraction using fixed fraction
      bisection refinement for the channel flow test problem.}
    \label{p2:fig:marking}
  \end{center}
\end{figure}

To study the effect of the choice of marking fraction, we plot in
Figure~\ref{p2:fig:marking} errors and efficiency indices for marking
fractions 0.1, 0.2, 0.3, 0.4 and 0.5 for fixed fraction bisection
refinement. We note that while a smaller marking fraction gives rise
to more efficient meshes, that is, a smaller number of degrees of
freedom are needed to reach a given level of accuracy, more refinement
levels are needed to reach that level of accuracy.

\subsubsection{Convergence of the global adaptive algorithm}

We next consider the convergence of the global adaptive algorithm. A
tolerance $\mathrm{TOL} = 0.001$ is prescribed for the error in the
goal functional, here the shear stress goal functional~$\M_1$ defined
in~\eqref{p2:eq:goal0}, and ask the global adaptive algorithm
described in Section~\ref{p2:sec:adaptive-algorithm} to adaptively
refine the mesh and select adaptive time steps such that $|\M_1(\ee)|
< \mathrm{TOL}$.

\begin{figure}[tbp!]
  \begin{center}
    \includegraphics[width=\plotwidth]{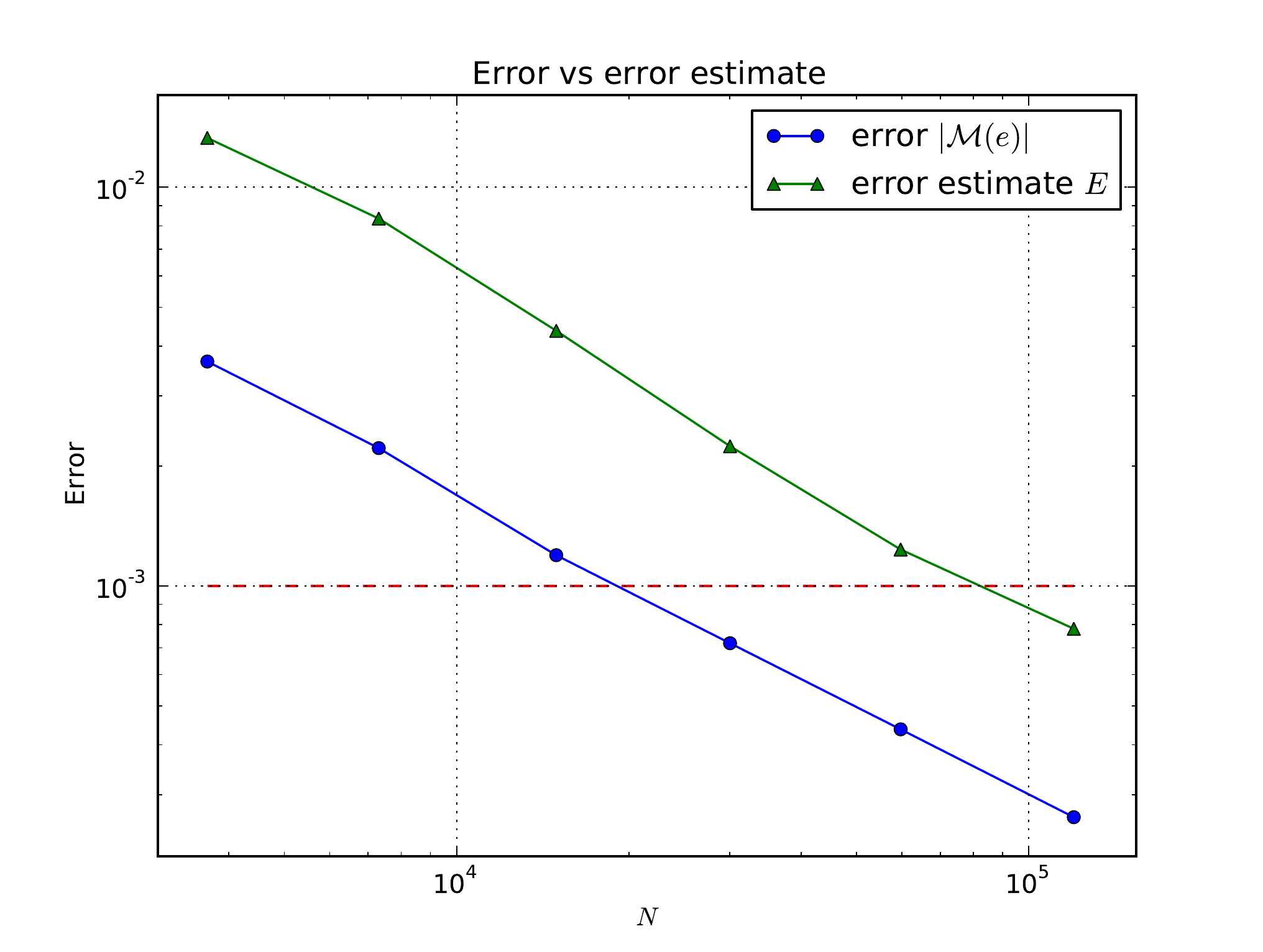}
    \includegraphics[width=\plotwidth]{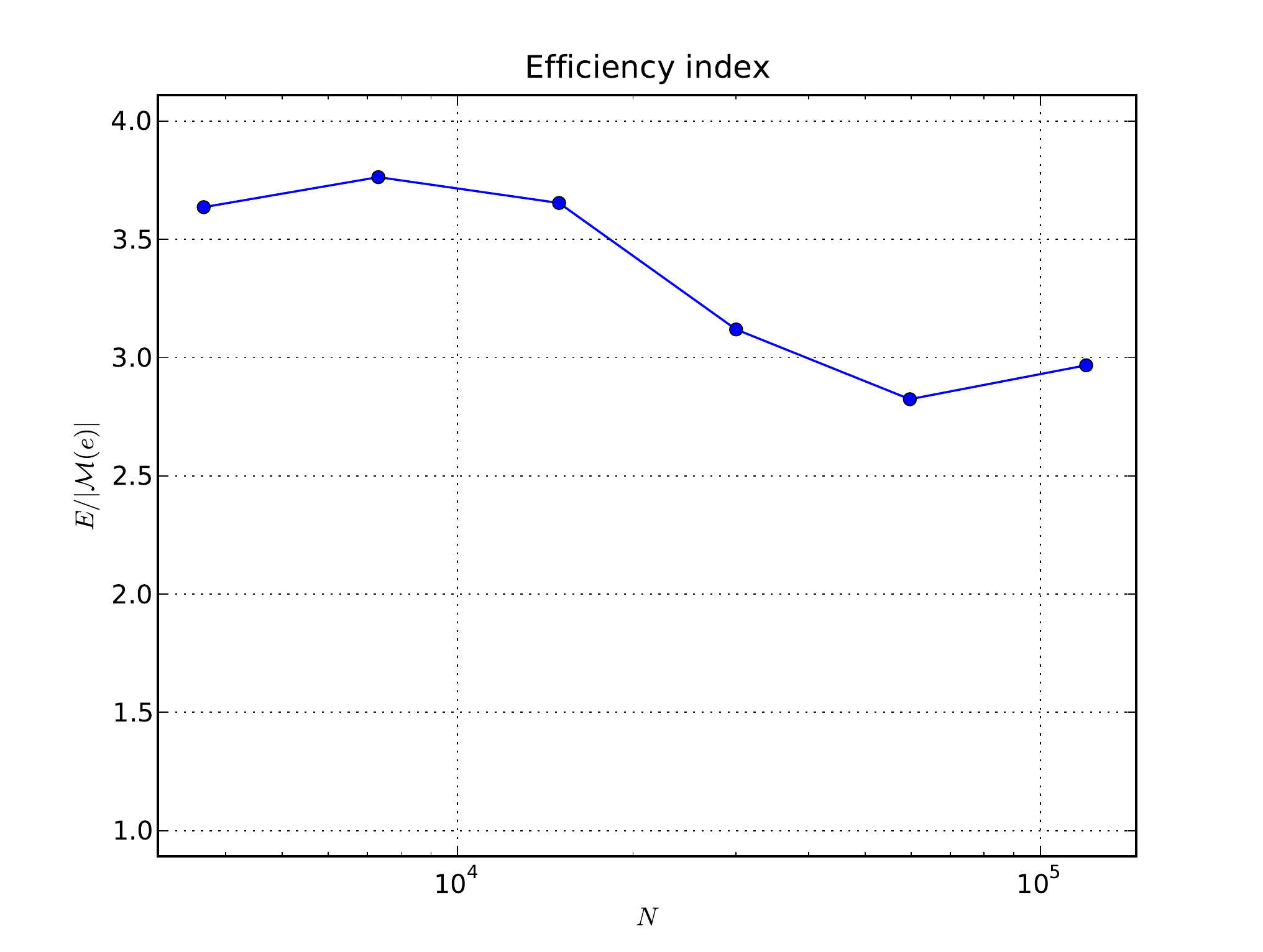}
    \caption{Convergence of the global space--time adaptive algorithm
      showing errors (top) and efficiency indices (bottom) using
      regular cut refinement with marking fraction $0.3$ for the
      channel flow test problem. The given tolerance $\mathrm{TOL} =
      0.001$ is reached after five refinements.}
    \label{p2:fig:global,ei}
  \end{center}
\end{figure}

Figure~\ref{p2:fig:global,ei} shows the convergence of the global
adaptive algorithm. The algorithm converges in five iterations when
the prescribed tolerance of $\mathrm{TOL}=0.001$ has been
reached. Although the actual error reaches the prescribed tolerance
after only three refinements, the adaptive algorithm performs well;
the size of the efficiency index is ca.~3. The adaptive time steps are
shown in Figure~\ref{p2:fig:global,k}. At $t = 0$, the time step is
set to the smallest time step from the previous refinement
level. Since the solution is initially at rest, the time residual is
initially small which leads to an increase in the size of the time
steps. As the fluid is accelerated by the pressure gradient, the time
residual increases and the time step is consequently reduced.

 \begin{figure}[tbp!]
  \begin{center}
    \includegraphics[width=0.8\textwidth]{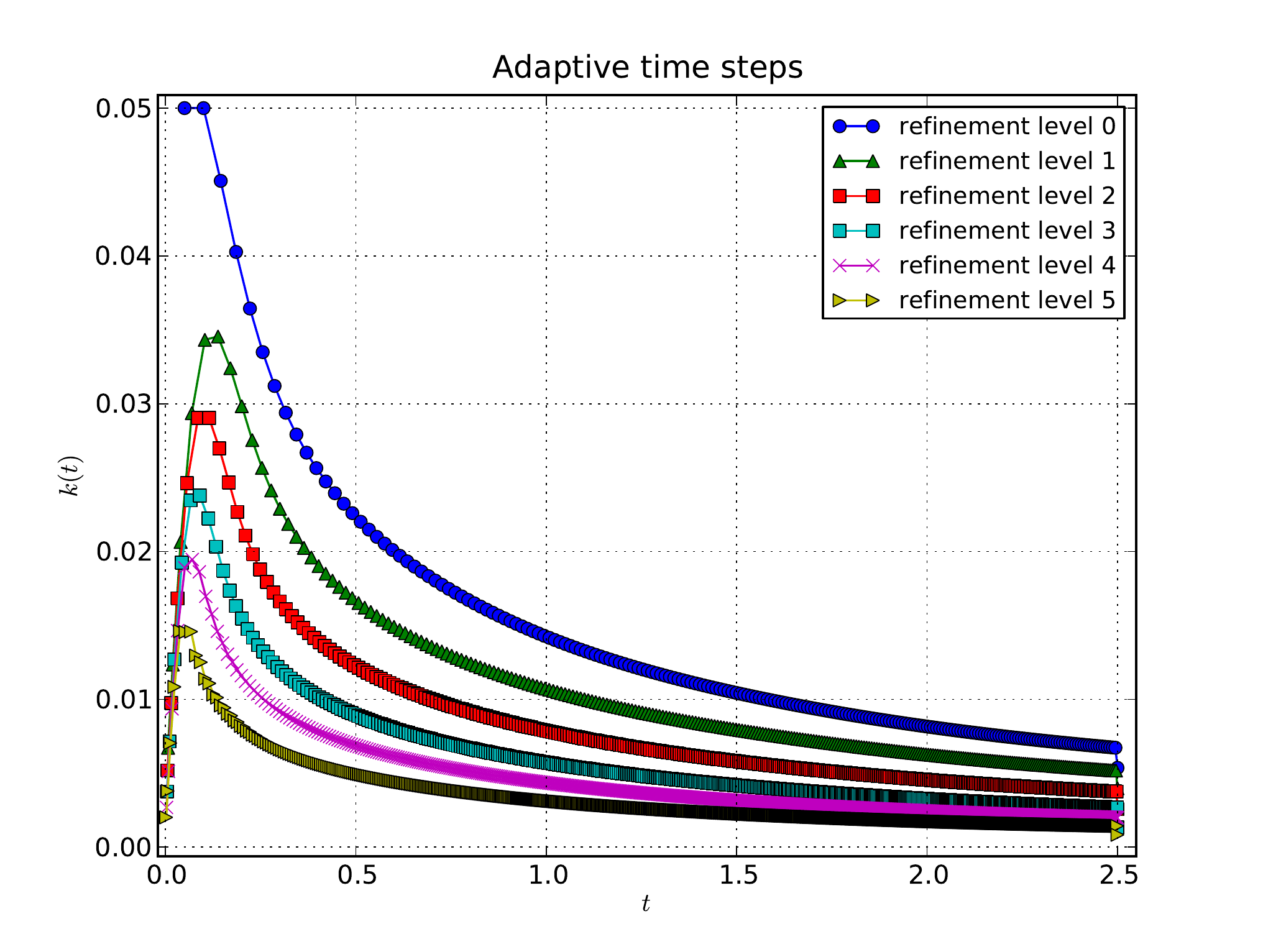}
    \vspace{-0.6cm}
    \caption{Time steps used by the global space--time adaptive
      algorithm on the time interval $[0, 2.5]$ for the channel flow
      test problem.}
    \label{p2:fig:global,k}
  \end{center}
\end{figure}

In Figures~\ref{p2:fig:global,contributions}
and~\ref{p2:fig:global,ec}, we plot the different contributions to the
total error estimate~$E = E_h + E_k + E_c$. We find that the error is
dominated by the space discretization error~$E_h$, while the time
discretization error remains small. This indicates that the time steps
are unnecessarily small. However, the time steps must remain small to
preserve stability of the numerical scheme. Although we have not taken
any special measures to control the size of the time step to maintain
numerical stability during mesh refinement, the adaptive time step
selection seems to decrease naturally in each adaptive iteration as a
result of an increase in the size of the residual $\dnorm{R}$. The
computational error $E_c$ remains practically constant throughout the
refinement and we note from Figure~\ref{p2:fig:global,ec} that the
dominating contribution to the computational error is from the
momentum equation; the discrete residual of the continuity equation
remains small.

\begin{figure}[tbp!]
  \begin{center}
    \includegraphics[width=0.8\textwidth]{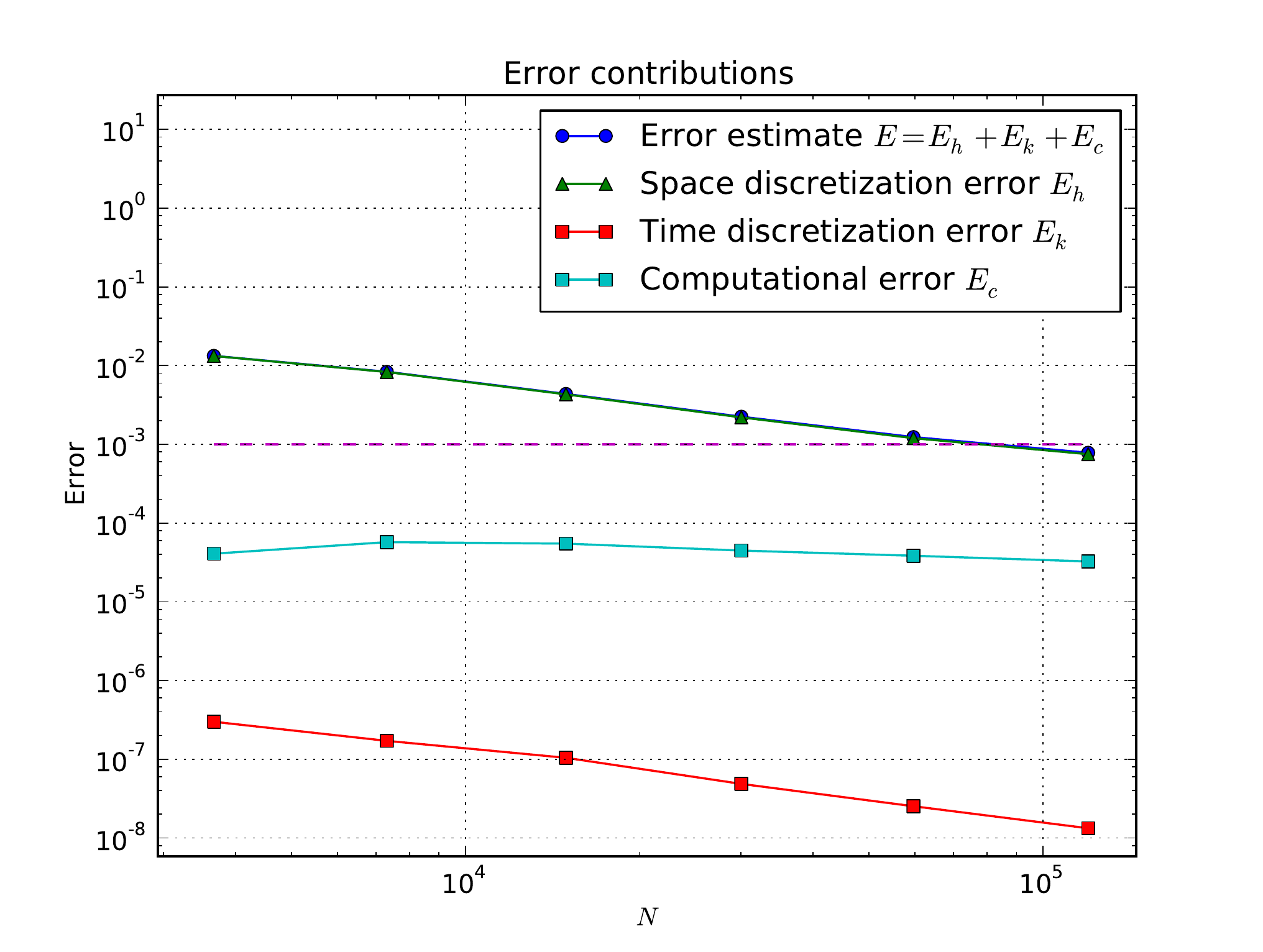}
    \caption{Contributions to the total error $E$ from spatial
      discretization ($E_h$), time discretization ($E_k$) and
      computational (splitting) error ($E_c$) for the global
      space--time adaptive algorithm using regular cut refinement with
      marking fraction $0.3$ for the channel flow test problem.}
    \label{p2:fig:global,contributions}
  \end{center}
\end{figure}

\begin{figure}[tbp!]
  \begin{center}
    \includegraphics[width=\plotwidth]{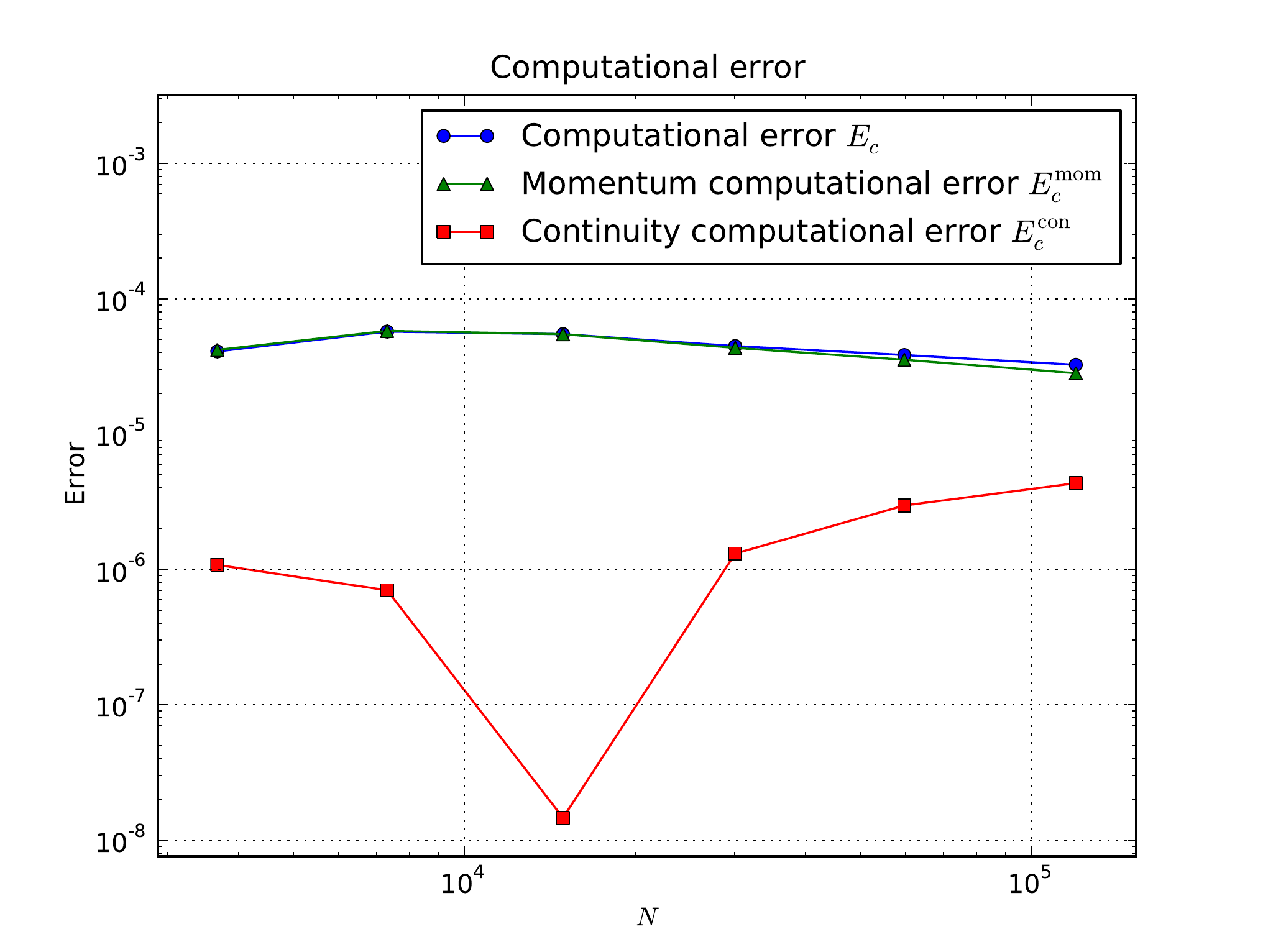}
    \caption{Contributions to the total computational (splitting)
      error $E_c$ from inexact solution of the finite element
      formulation of the momentum equation ($E_c^{\mathrm{mom}}$) and
      the continuity equation ($E_c^{\mathrm{con}}$) using regular cut
      refinement with marking fraction $0.3$ for the channel flow test
      problem.}
    \label{p2:fig:global,ec}
  \end{center}
\end{figure}

\subsubsection{Convergence as function of $h$ and $k$}

Finally, we investigate how the error contributions $E_h$, $E_k$ and
$E_c$ depend on the mesh size $h$ and the time step $k$. We consider
the shear stress goal functional~$\M_1$ defined in~\eqref{p2:eq:goal0}
computed on a sequence of uniformly refined meshes with mesh sizes $h
= 0.2$, $h = 0.1$, $h = 0.05$ and $h = 0.025$, and fixed time steps $k
= 0.01$, $k = 0.005$, $k = 0.0025$ and $k = 0.00125$.

Figure~\ref{p2:fig:hk,Eh} shows the space discretization error $E_h$
as function of the mesh size~$h$. The results indicate that the
convergence of the error in the goal functional is linear with respect
to the mesh size. This has not been considered in detail but we note
that for a $P_2$--$P_1$ Taylor--Hood discretization, we expect the
convergence of the error in the velocity to be $h^3$ in the mesh
size. However, as the goal functional $\M_1$ involves the shear
stress, the order of convergence is reduced to $h^2$. The convergence
rate is further decreased by the fact that the goal functional
considers the shear stress on the boundary and as a result of the
singularities at the reentrant corners close to the evaluation of the
goal functional. We further note from this figure that $E_h$ does not
depend on the size of the time step with one exception; the error goes
up on the finest mesh for the largest time step $k = 0.01$, indicating
instability of the numerical scheme for large relative time steps.

In Figure~\ref{p2:fig:hk,Ek}, we plot the time discretization error
$E_k$ as function of mesh size $h$ and time step $k$, respectively.
We conclude that $E_k$ depends only weakly on $h$ and that the
convergence of $E_k$ is quadratic in the time step $k$.

For the computational error $E_c$ displayed in
Figures~\ref{p2:fig:hk,Ec0} and~\ref{p2:fig:hk,Ec1}, we similarly find
a weak dependence on the mesh size $h$. We further note that the
contribution from the momentum equation is linear in the size of the
time step, whereas the contribution from the continuity equation is
quadratic. Overall, we thus find that the order of convergence is
linear in the time step as expected.

\begin{figure}
  \begin{center}
    \includegraphics[width=\plotwidth]{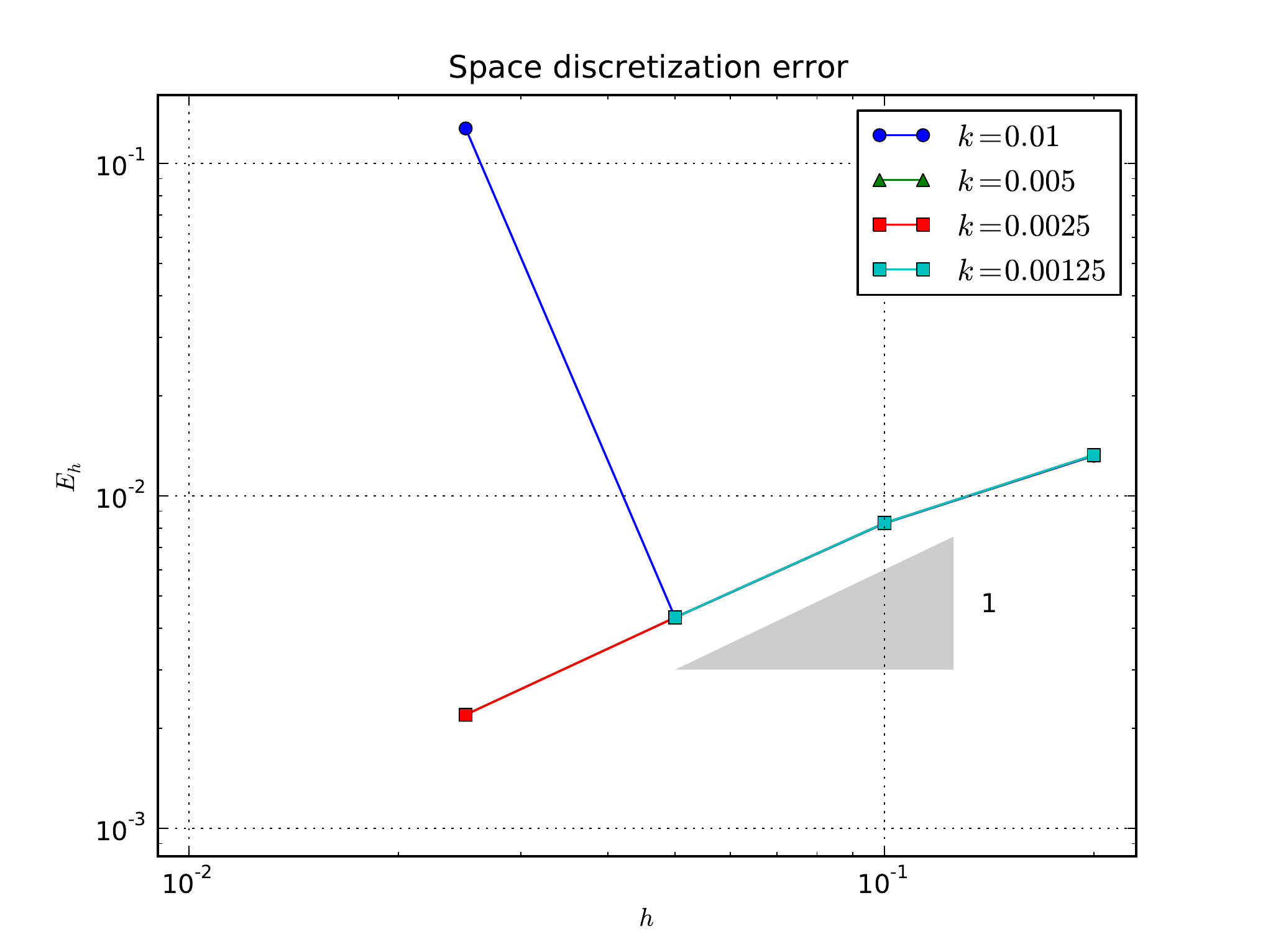}
    \caption{Space discretization error $E_h$ for the shear stress
      goal functional~$\M_1$ defined in~\eqref{p2:eq:goal0} as
      function of mesh size~$h$ for varying (fixed) time steps.}
    \label{p2:fig:hk,Eh}
  \end{center}
\end{figure}

\begin{figure}
  \begin{center}
    \includegraphics[width=\plotwidth]{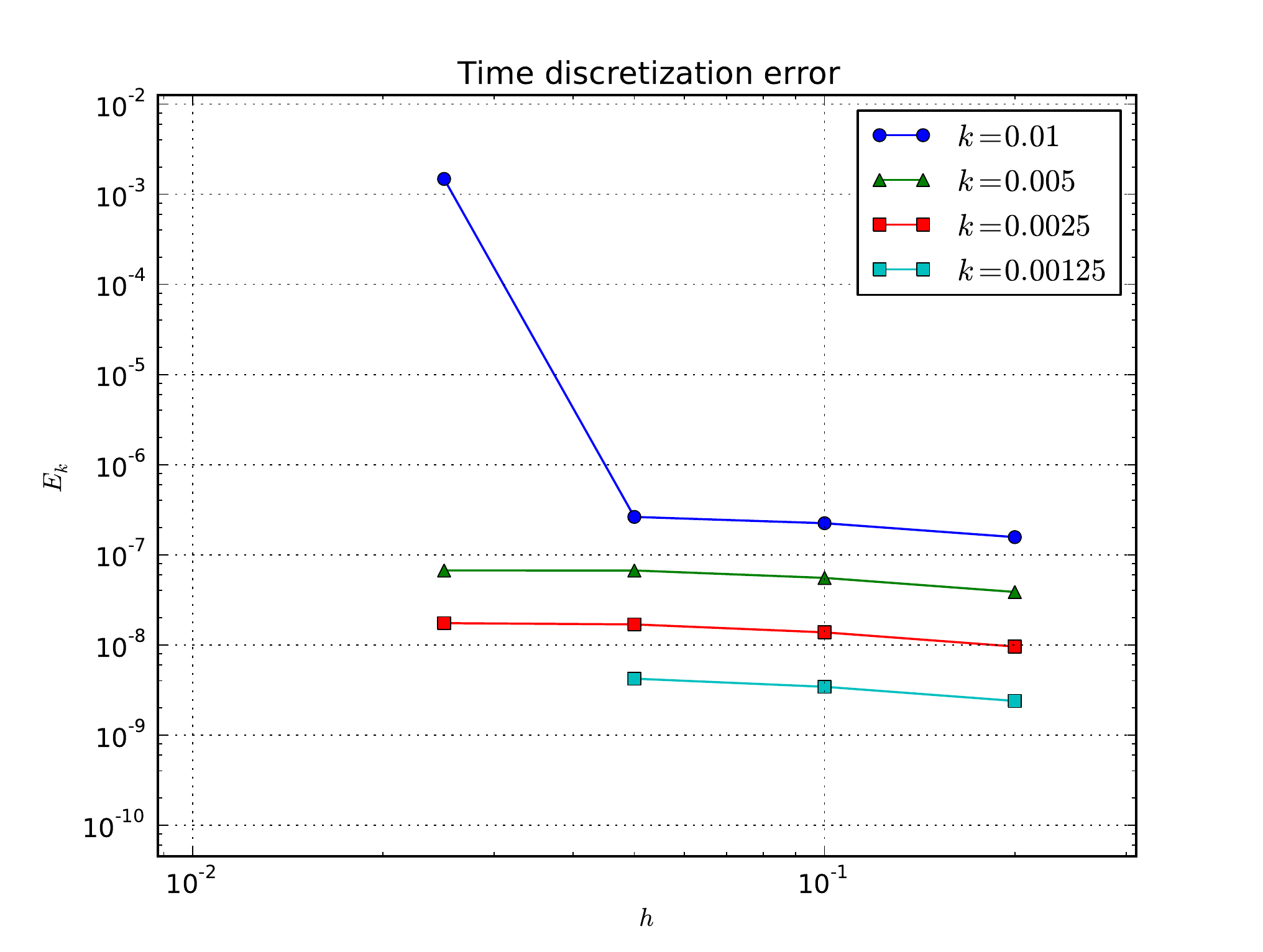}
    \includegraphics[width=\plotwidth]{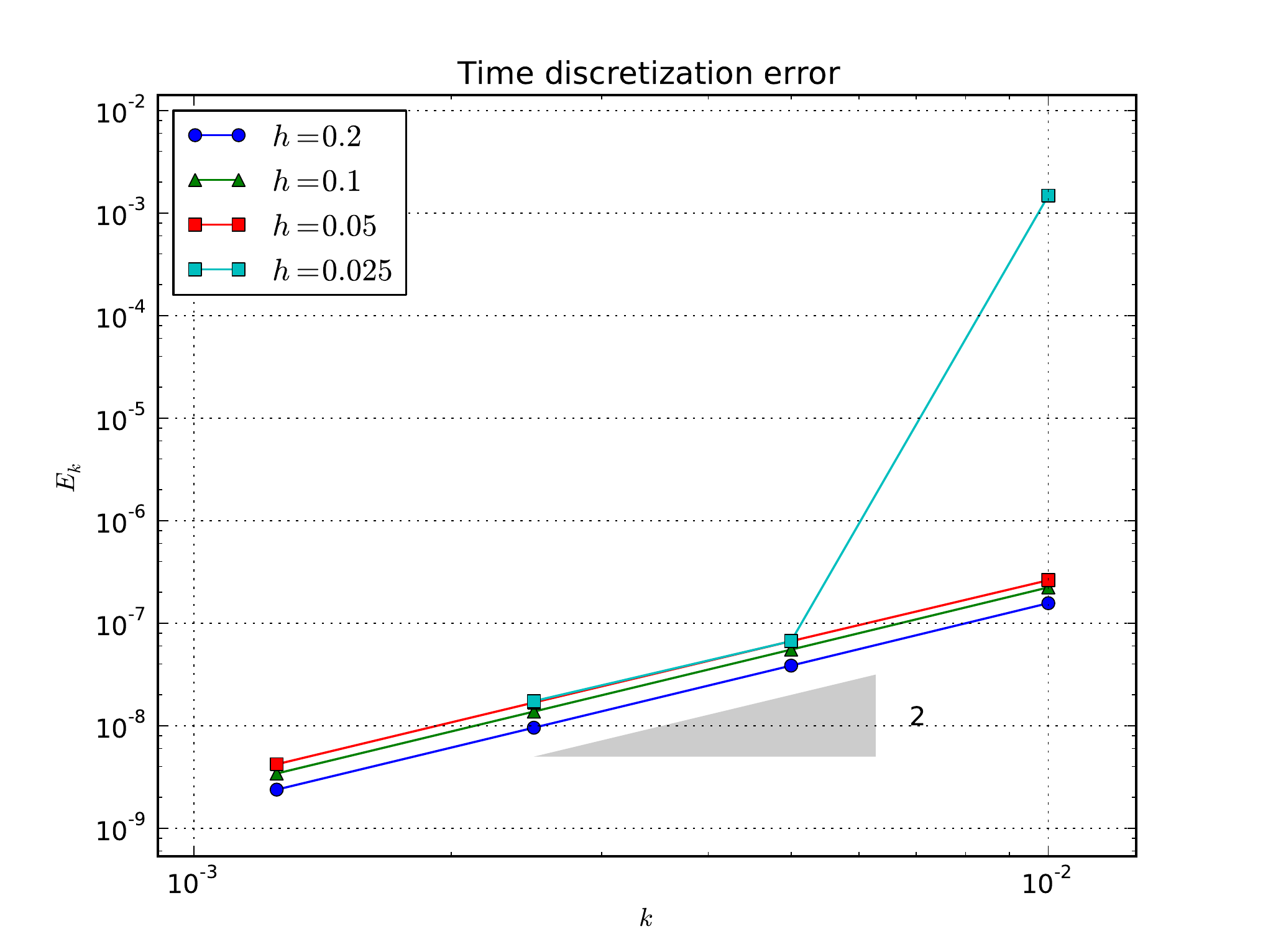}
    \caption{Time discretization error $E_k$ for the shear stress goal
      functional~$\M_1$ defined in~\eqref{p2:eq:goal0} as function of
      mesh size~$h$ (top) and time step size~$k$ (bottom).}
    \label{p2:fig:hk,Ek}
  \end{center}
\end{figure}

\begin{figure}
  \begin{center}
    \includegraphics[width=\plotwidth]{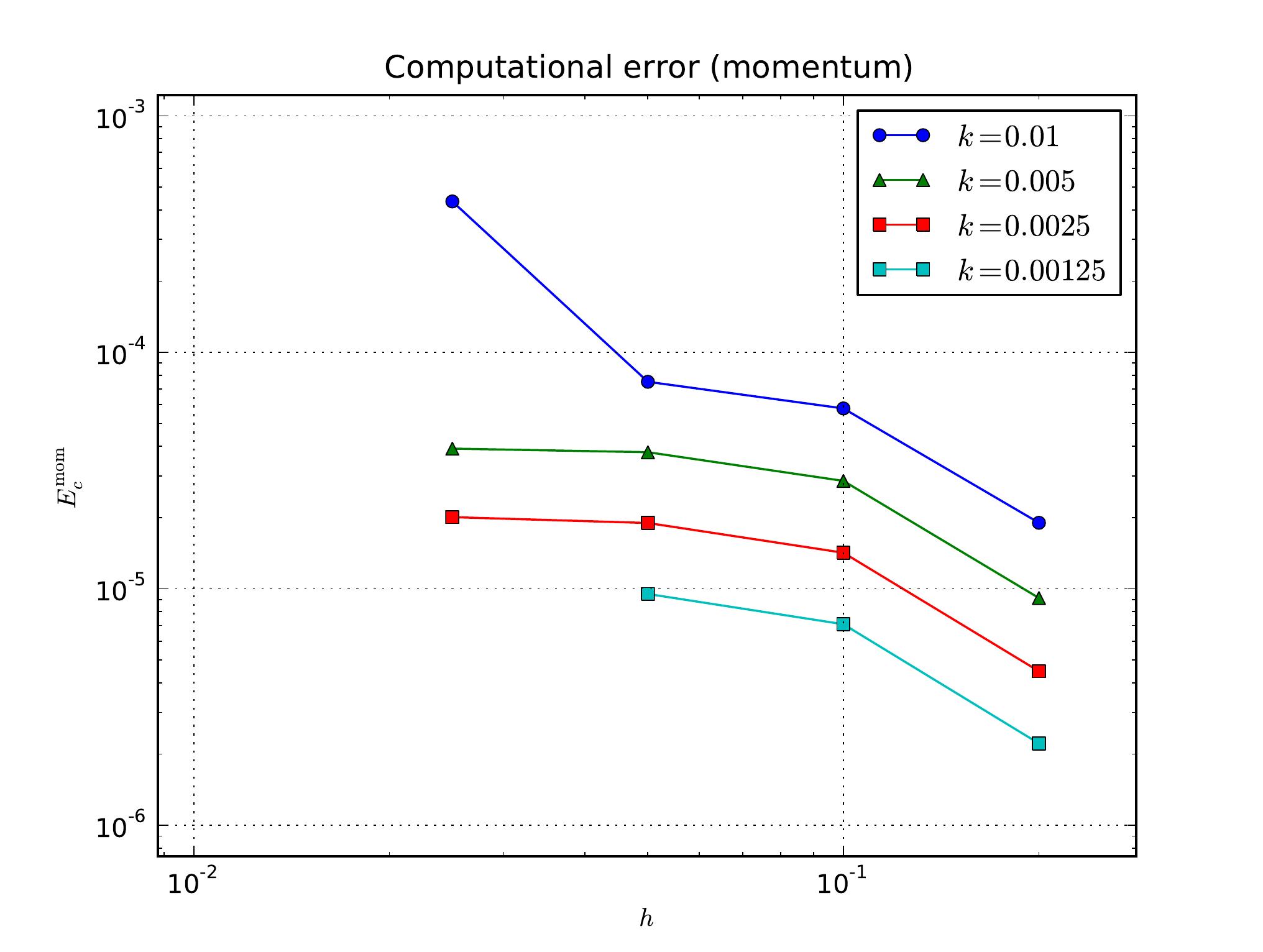}
    \includegraphics[width=\plotwidth]{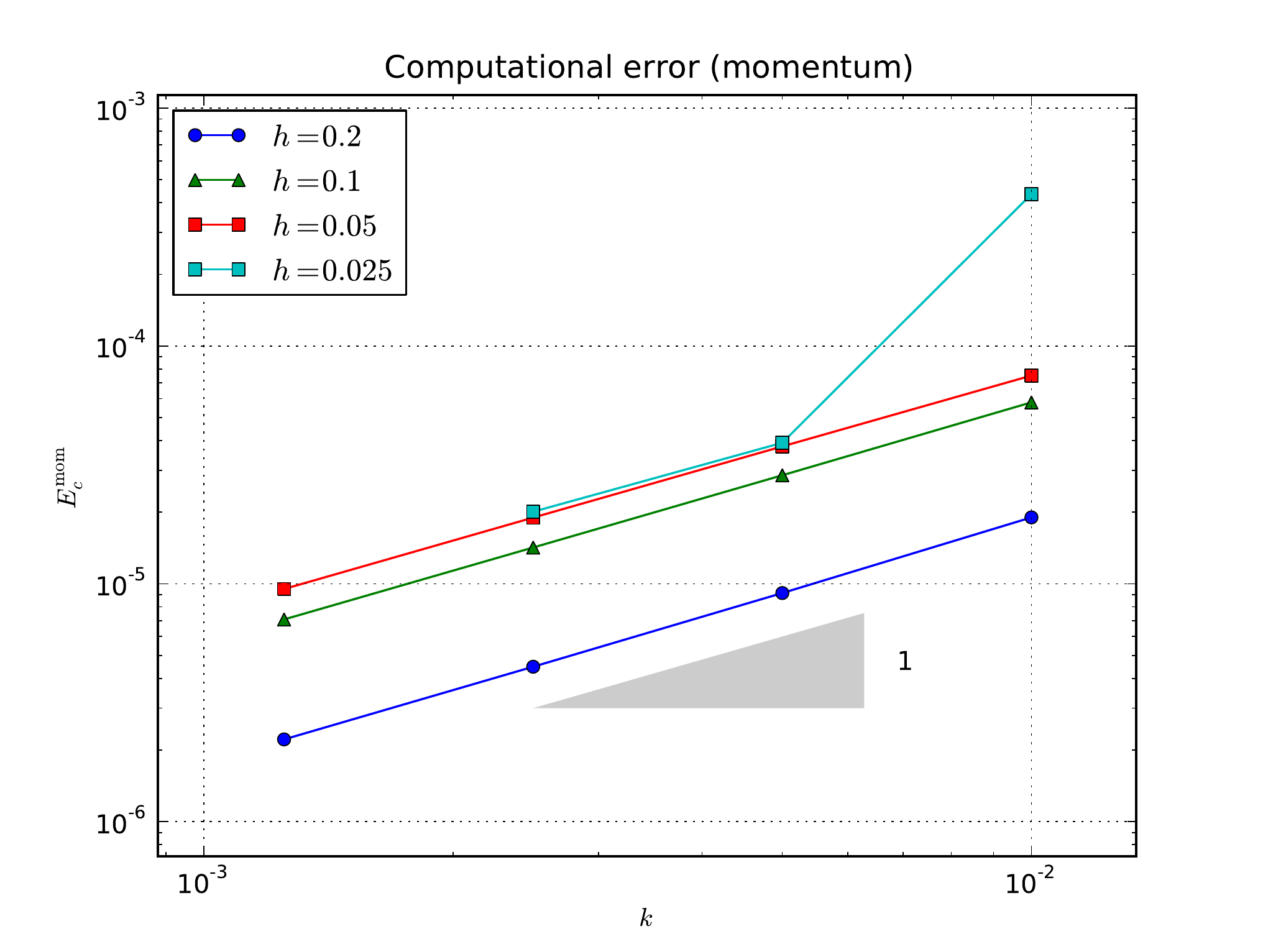}
    \caption{Momentum computational error $E_c^{\mathrm{mom}}$ for the
      shear stress goal functional~$\M_1$ defined
      in~\eqref{p2:eq:goal0} as function of mesh size~$h$ (top) and
      time step size~$k$ (bottom).}
    \label{p2:fig:hk,Ec0}
  \end{center}
\end{figure}

\begin{figure}
  \begin{center}
    \includegraphics[width=\plotwidth]{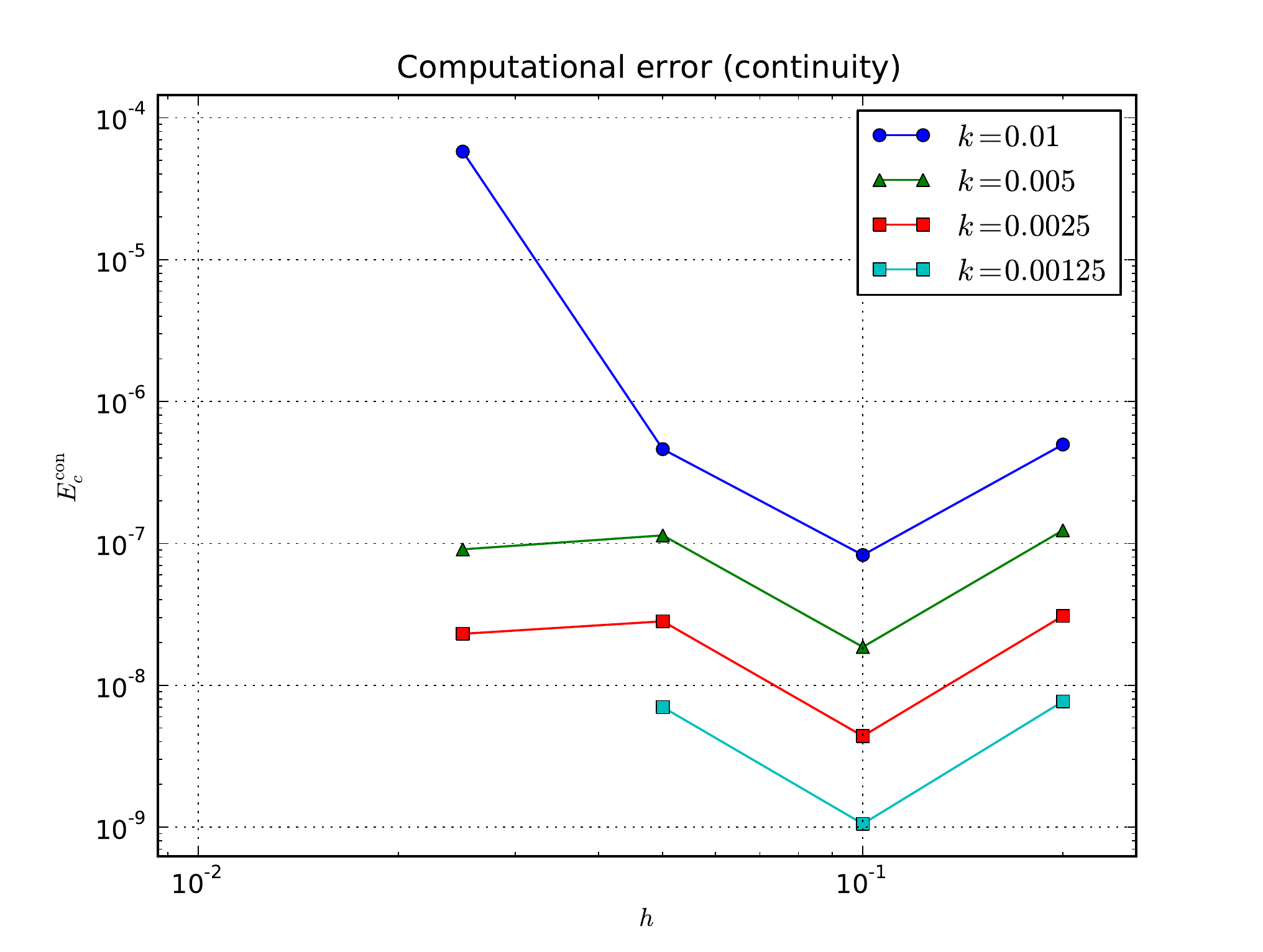}
    \includegraphics[width=\plotwidth]{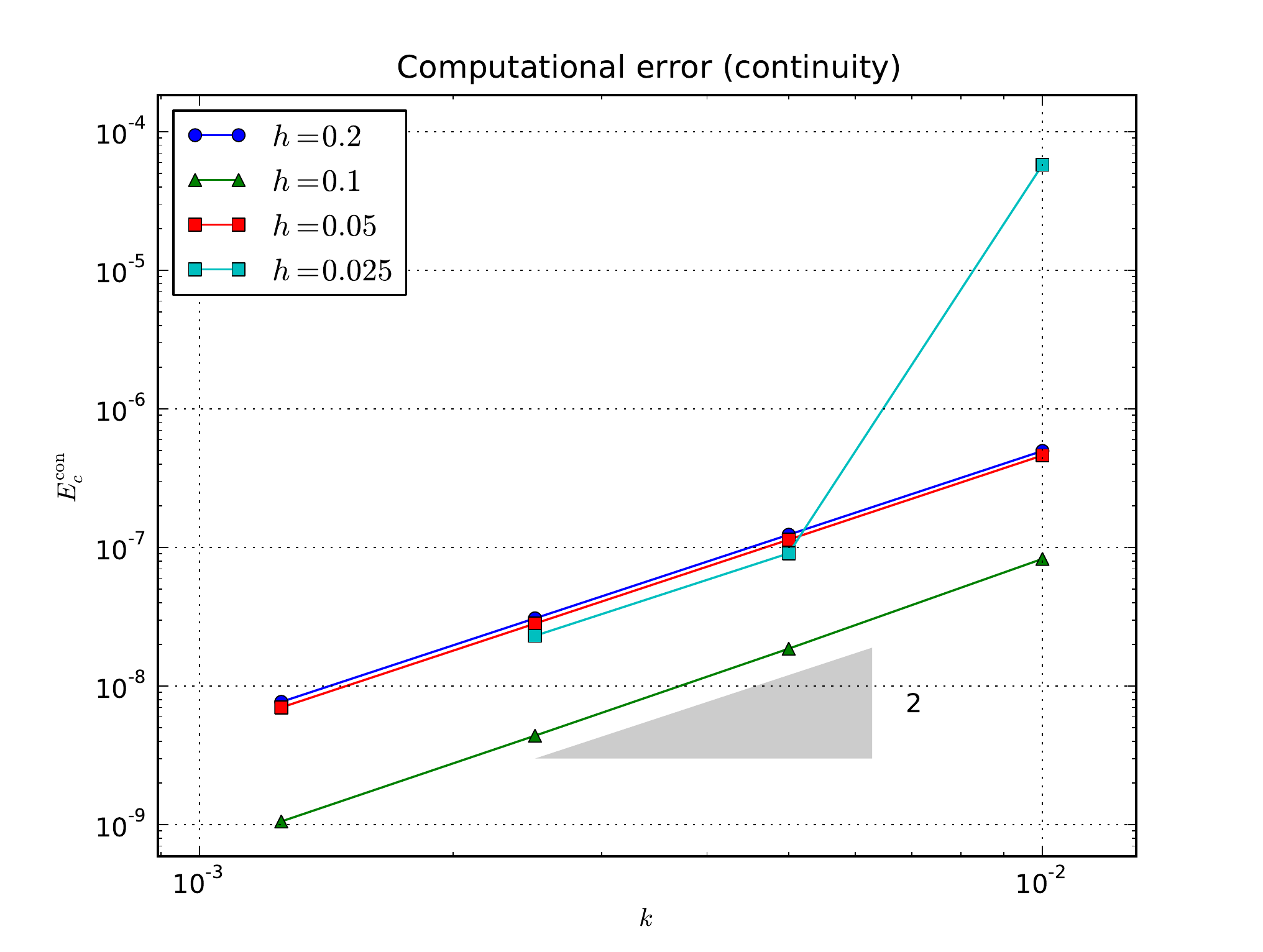}
    \caption{Continuity computational error $E_c^{\mathrm{con}}$ for
      the shear stress goal functional~$\M_1$ defined
      in~\eqref{p2:eq:goal0} as function of mesh size~$h$ (top) and
      time step size~$k$ (bottom).}
    \label{p2:fig:hk,Ec1}
  \end{center}
\end{figure}

\subsection{Case 2: Lid-driven cavity}

As a second test problem, we consider the lid-driven cavity problem on
the unit square $(0, 1) \times (0, 1)$. As boundary conditions, we set
$u = (x_1 (1 - x_1), 0)$ at the top of the cavity ($x_2 = 1$) with no-slip
boundary conditions on the remaining boundary for the velocity. We
also fix the pressure $p = 0$ at the bottom of the cavity ($x_2 = 0$).
This ``unphysical'' boundary condition for the pressure gives rise to
(small) gradients in the pressure field in the vicinity of $x_2 = 0$. A
better way to ensure solvability of the pressure update step of
Algorithm~\ref{p2:alg:ipcs} is to require $\int_{\Omega} p \dx =
0$. However, we have here chosen to use a Dirichlet boundary condition
for the pressure, as this is often used in applications and we wish to
study its effect on mesh refinement.

We set the kinematic viscosity to $\nu = 1$ and run the simulation
over the time interval $[0, 1]$. As a goal functional, we consider a
Gaussian-weighted average of the $x_2$-component of the velocity field:
\begin{equation} \label{p2:eq:goal3}
  \M_2(\uu)
  = \int_0^T \int_{\Omega} u_2(x, t) \, \phi(x) \dx \dt.
\end{equation}
The weight function $\phi$ is chosen as
\begin{equation}
  \phi(x_1, x_2) = c \exp(-((x_1 - \bar{x}_1)^2 + (x_2 - \bar{x}_2)^2) / 2r^2),
\end{equation}
where $(\bar{x}_1, \bar{x}_2) = (0.75, 0.75)$, $r = 0.15$ and $c
\approx 27.571034$ is chosen such that $\int_{\Omega} \phi(x_1, x_2)
\dx = 1$. As a reference value, we take $\M_2(\uu) = -0.039389$. The
velocity and pressure fields at final time $T = 1$ are shown in
Figure~\ref{p2:fig:driven_cavity,solution}.

\begin{figure}[tbp!]
  \begin{center}
    \includegraphics[width=5.9cm]{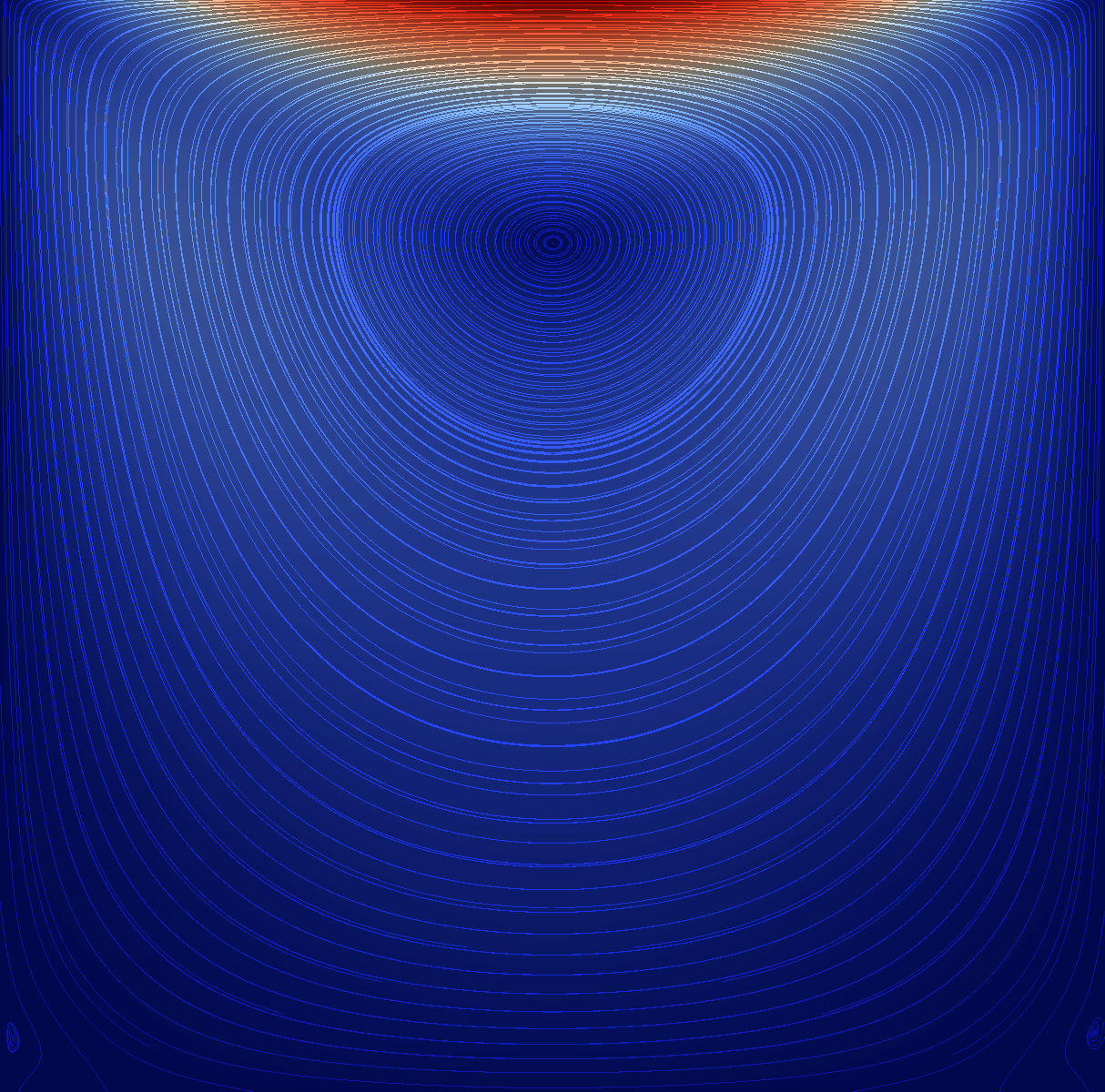}
    \includegraphics[width=5.9cm]{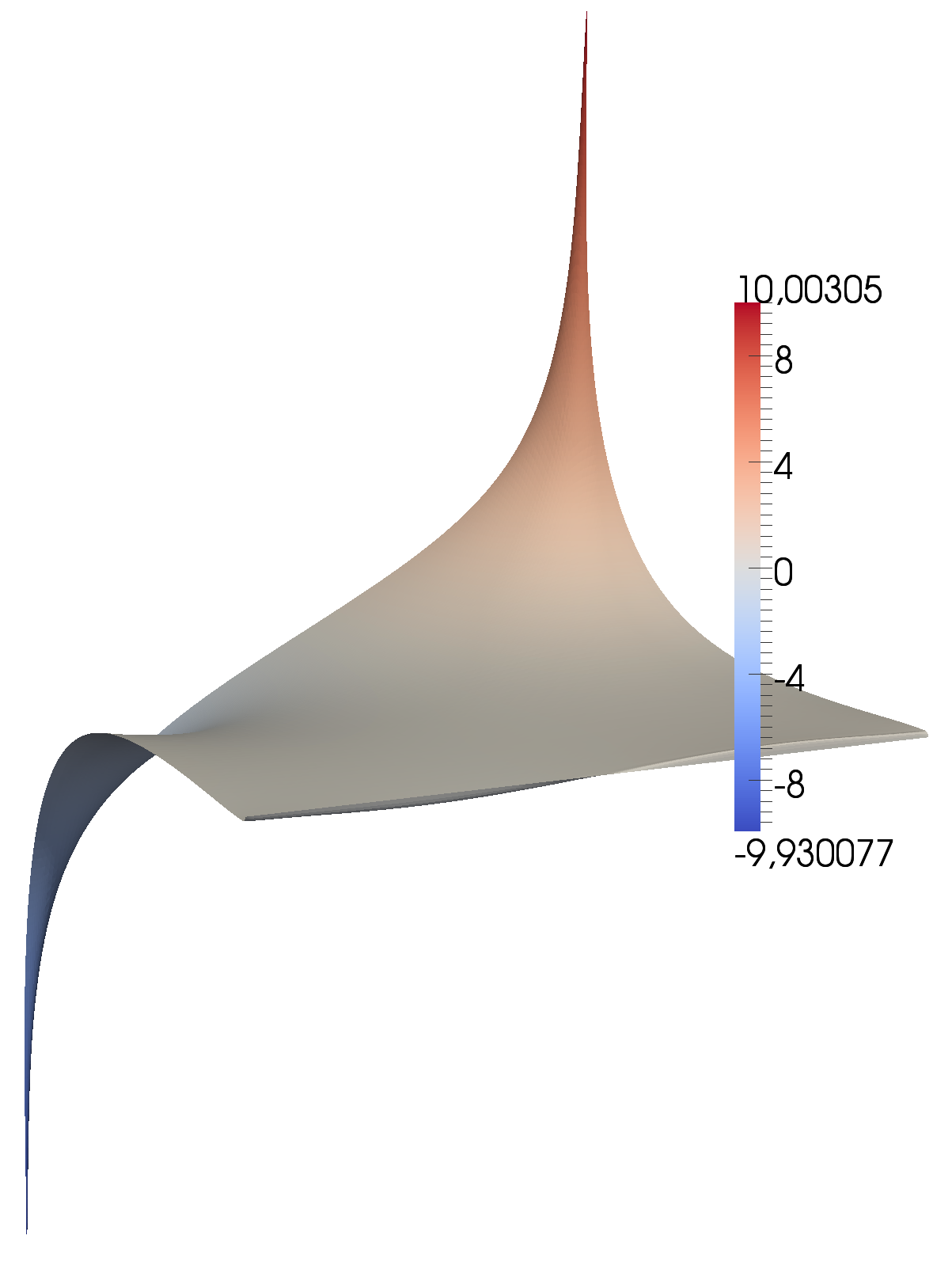}
    \caption{Velocity field (left) and pressure (right) at final time
      $T = 1$ for the lid-driven cavity test problem.}
    \label{p2:fig:driven_cavity,solution}
  \end{center}
\end{figure}

\subsubsection{Dual solutions and adaptive meshes}

The choice of goal functional generates a source located in the top
right corner (at $x_1 = x_2 = 0.75$). The dual solution is advected
backwards along the primal velocity field and the resulting dual
velocity field is shown in
Figure~\ref{p2:fig:driven_cavity,dual}. Notice the large secondary
vortex in the top right corner and the small secondary vortices in the
other three corners. The corresponding adaptive mesh is refined
heavily in the top left and right corners (see
Figure~\ref{p2:fig:driven_cavity,mesh}), as a result of pressure
spikes in these corners, but also at the bottom of the cavity as a
result of the Dirichlet boundary condition used for the pressure.

\begin{figure}[tbp!]
  \begin{center}
    \includegraphics[width=7cm]{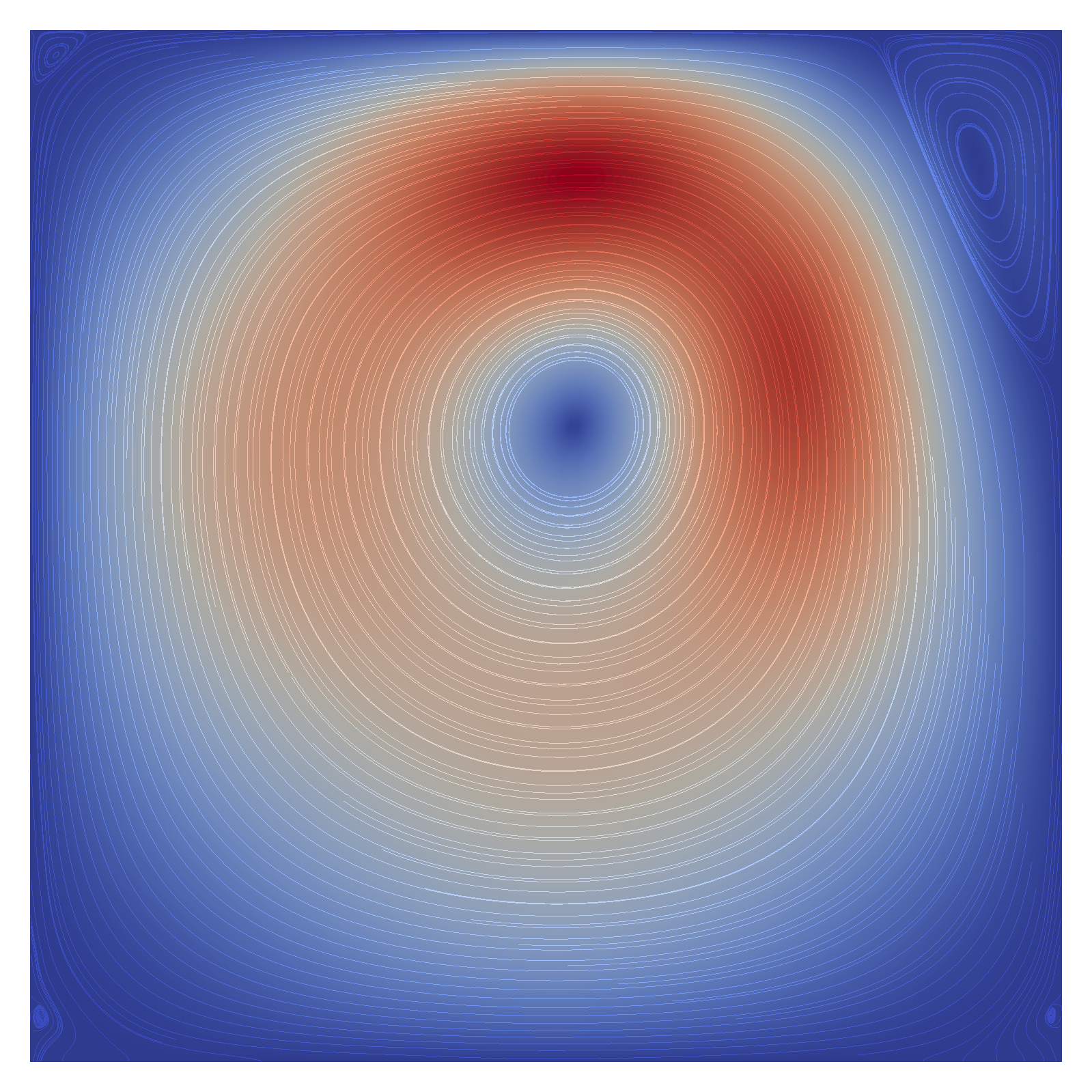}
    \caption{Dual fluid velocity field at ``final'' time $t = 0$ for
      the lid-driven cavity test problem.}
    \label{p2:fig:driven_cavity,dual}
  \end{center}
\end{figure}

\begin{figure}[tbp!]
  \begin{center}
    \includegraphics[width=5.9cm]{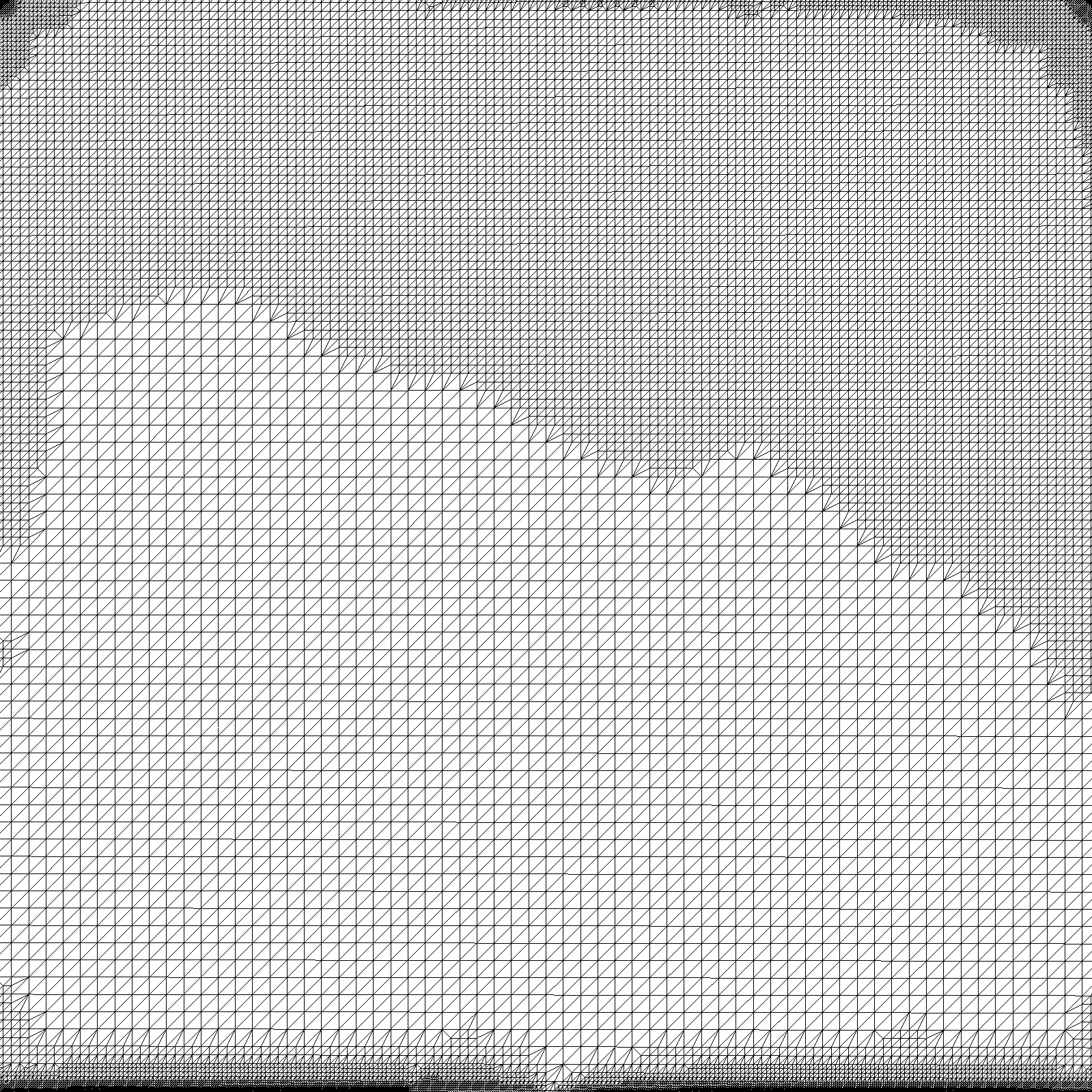}
    \includegraphics[width=5.9cm]{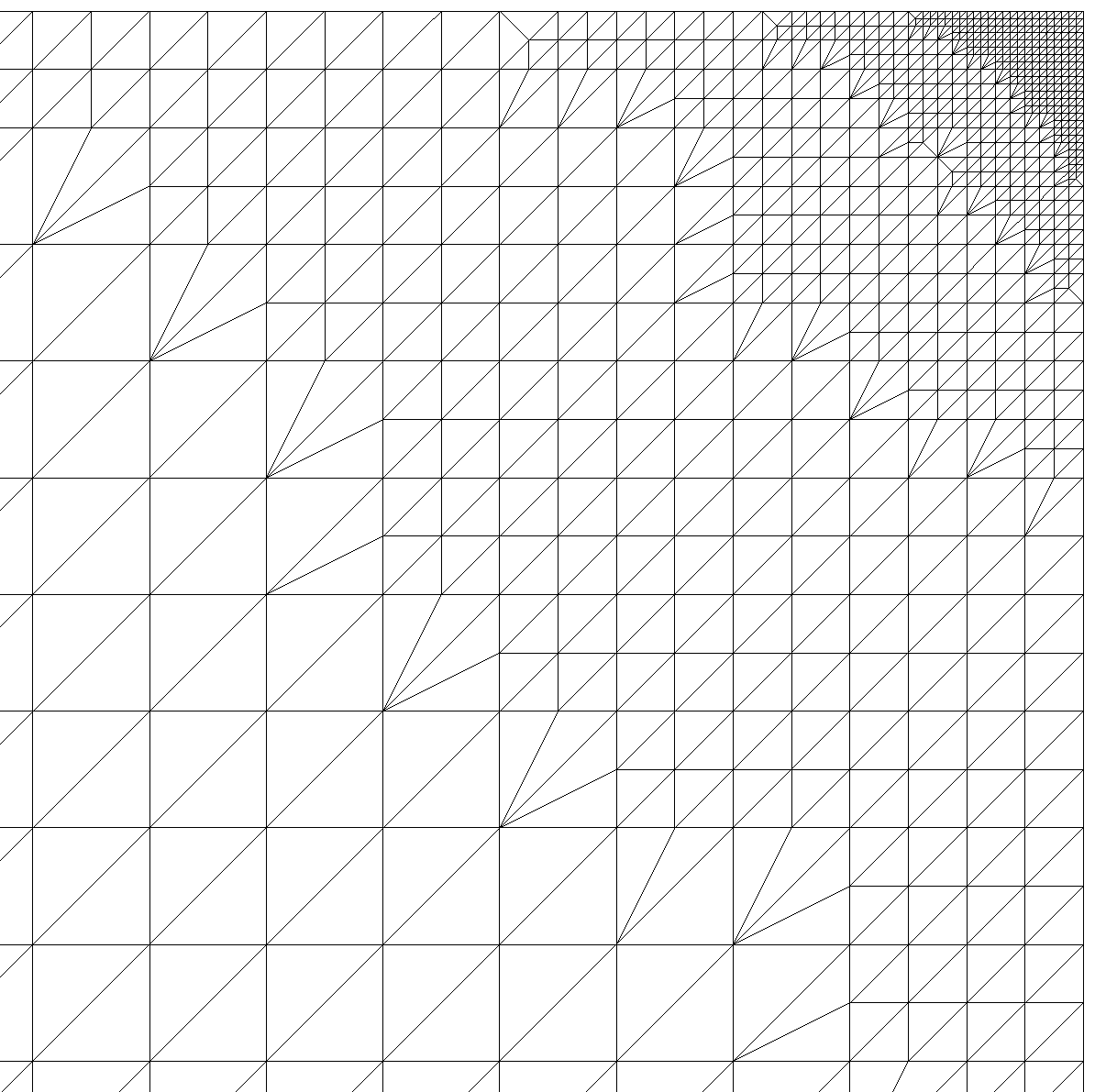}
    \caption{Mesh obtained by 12 levels of regular cut refinement with
      marking fraction 0.3 for the lid-driven cavity test problem
      (left) and a detailed view of the refined mesh in the top right
      corner (right).}
    \label{p2:fig:driven_cavity,mesh}
  \end{center}
\end{figure}

\subsubsection{Error and efficiency indices}

Figure~\ref{p2:fig:driven_cavity,ei} shows the error of the goal
functional and the corresponding efficiency indices for a sequence of
adaptively refined meshes, using adaptive time-stepping on each
refined mesh. Two different adaptive refinement algorithms, recursive
bisection and regular cut refinement, are compared to uniform
refinement. Both adaptive algorithms perform significantly better
compared to uniform refinement. No significant difference can be noted
between the two adaptive refinement algorithms, other than that
recursive bisection requires approximately twice the number of
refinement levels to reach the same level of accuracy as regular cut
refinement. We further note that the efficiency indices vary between
ca.~$1$ and~$10$. Interestingly, the efficiency indices for uniform
refinement seem to converge towards~$1$.

\begin{figure}[tbp!]
  \begin{center}
    \includegraphics[width=\plotwidth]{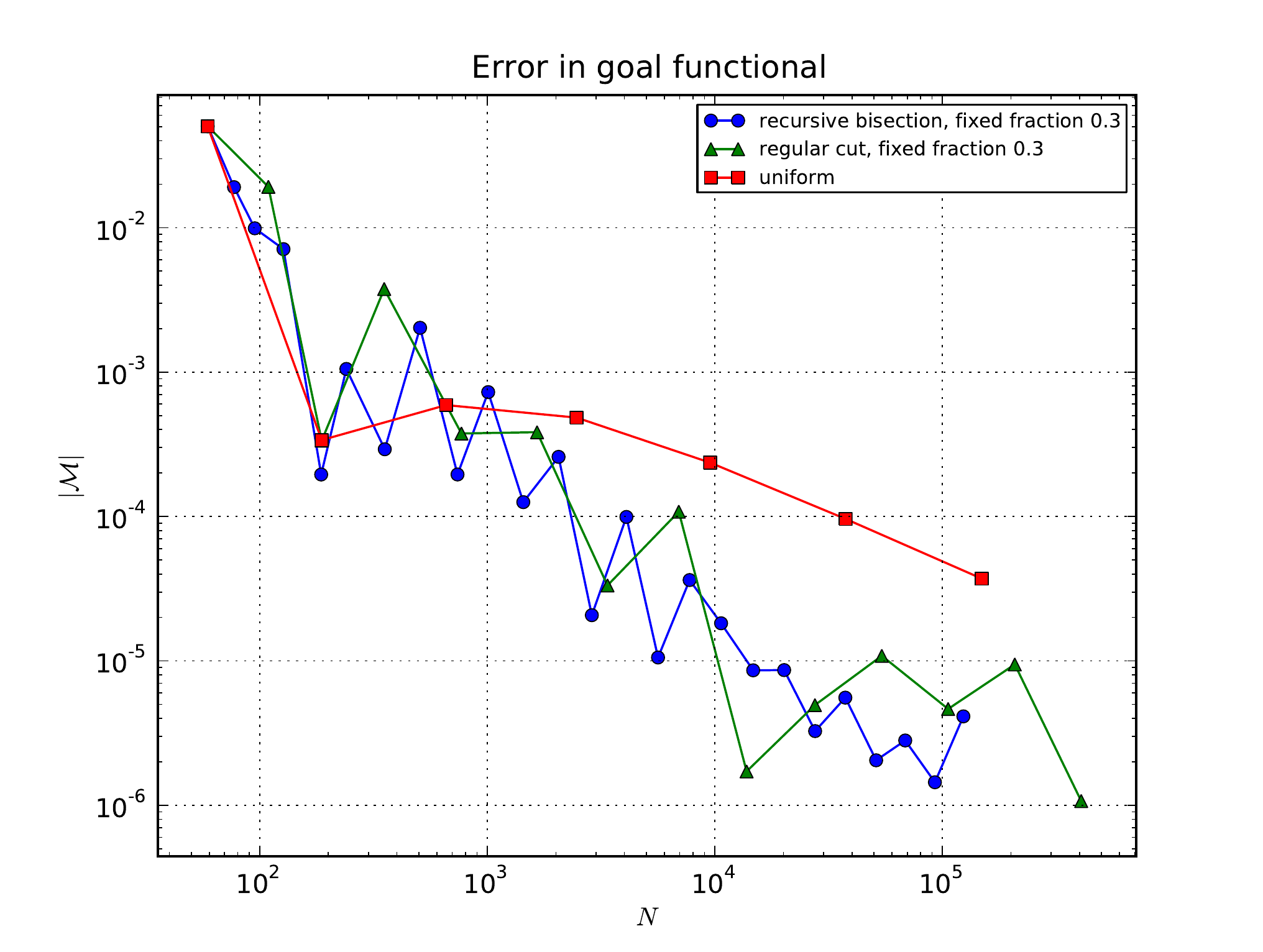}
    \includegraphics[width=\plotwidth]{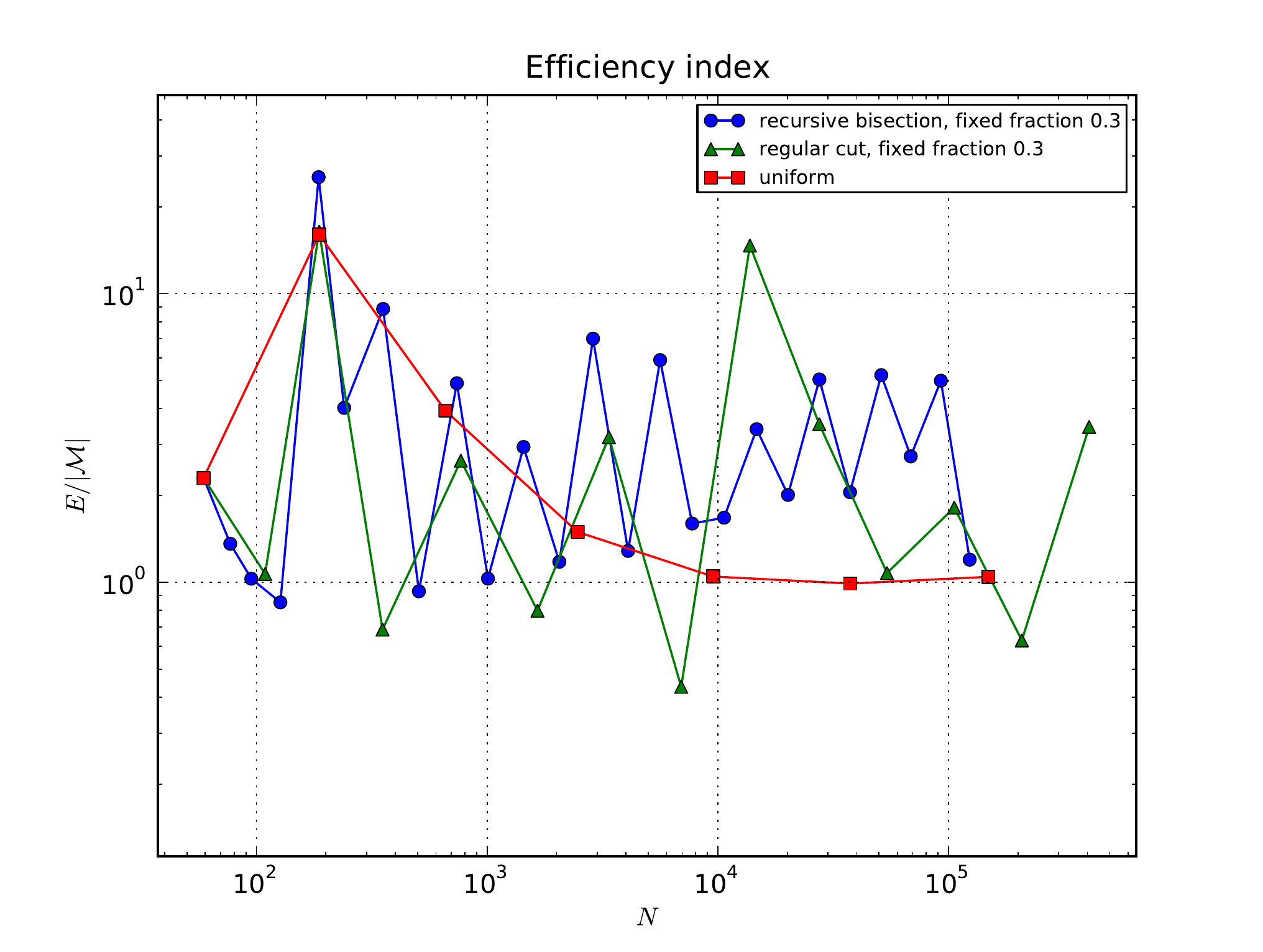}
    \caption{Error (top) and efficiency indices (bottom) as function
      of the number of spatial degrees of freedom using adaptive time
      steps, fixed fraction marking (marking fraction $0.3$) and
      varying refinement algorithms (recursive bisection, regular cut
      and uniform) for the lid-driven cavity problem.}
    \label{p2:fig:driven_cavity,ei}
  \end{center}
\end{figure}

\section{Conclusions}
\label{p2:sec:conclusion}

We have presented an adaptive finite element method for the
incompressible Navier--Stokes equations based on a standard splitting
scheme (incremental pressure correction). By treating the splitting
scheme as an approximation of a pure Galerkin finite element scheme,
one may analyze the error as a sum of contributions from space
discretization, time discretization and a computational error that
measures the deviation of the splitting scheme from the pure Galerkin
scheme. Numerical experiments indicate good performance of the
adaptive algorithm and error estimates that closely match the actual
error. The proposed method may thus serve as an attractive approach to
solving the incompressible Navier--Stokes equations, combining the
efficiency of a simple splitting method with the framework of
goal-oriented adaptive finite element methods.

The presented adaptive algorithm can be further improved by extending
the adaptive time step selection to control the size of the
computational error $E_c$. It may also be interesting to consider
modified splitting schemes to reduce the size of the computational
error, in particular the size of the discrete momentum residual.

\section{Acknowledgments}

This work is supported by an Outstanding Young Investigator grant from
the Research Council of Norway, NFR 180450. This work is also
supported by a Center of Excellence grant from the Research Council of
Norway to the Center for Biomedical Computing at Simula Research
Laboratory.

\bibliographystyle{model1-num-names}
\bibliography{bibliography}

\end{document}

%% file: pdf/channel_with_flap.pdf_tex

\begingroup
  \makeatletter
  \providecommand\color[2][]{%
    \errmessage{(Inkscape) Color is used for the text in Inkscape, but the package 'color.sty' is not loaded}
    \renewcommand\color[2][]{}%
  }
  \providecommand\transparent[1]{%
    \errmessage{(Inkscape) Transparency is used (non-zero) for the text in Inkscape, but the package 'transparent.sty' is not loaded}
    \renewcommand\transparent[1]{}%
  }
  \providecommand\rotatebox[2]{#2}
  \ifx\svgwidth\undefined
    \setlength{\unitlength}{250.91469727pt}
  \else
    \setlength{\unitlength}{\svgwidth}
  \fi
  \global\let\svgwidth\undefined
  \makeatother
  \begin{picture}(1,0.25735215)%
    \put(0,0){\includegraphics[width=\unitlength]{channel_with_flap.pdf}}%
    \put(0.47073705,0.03197987){\color[rgb]{0,0,0}\makebox(0,0)[lb]{\smash{$4.0$}}}%
    \put(0.17422193,0.00328486){\color[rgb]{0,0,0}\makebox(0,0)[lb]{\smash{$1.4$}}}%
    \put(0.31769699,0.00328486){\color[rgb]{0,0,0}\makebox(0,0)[lb]{\smash{$0.4$}}}%
    \put(-0.00113647,0.16907825){\color[rgb]{0,0,0}\makebox(0,0)[lb]{\smash{$1.0$}}}%
    \put(0.38465201,0.13719491){\color[rgb]{0,0,0}\makebox(0,0)[lb]{\smash{$0.6$}}}%
    \put(0.05626543,0.18779104){\color[rgb]{0,0,0}\makebox(0,0)[lb]{\smash{$p=1$}}}%
    \put(0.69393234,0.18779104){\color[rgb]{0,0,0}\makebox(0,0)[lb]{\smash{$p=0$}}}%
    \put(0.05626543,0.13996602){\color[rgb]{0,0,0}\makebox(0,0)[lb]{\smash{$\nabla u \cdot n = 0$}}}%
    \put(0.69393234,0.13996602){\color[rgb]{0,0,0}\makebox(0,0)[lb]{\smash{$\nabla u \cdot n = 0$}}}%
    \put(0.5058206,0.23242772){\color[rgb]{0,0,0}\makebox(0,0)[lb]{\smash{$u = 0$}}}%
    \put(0.5058206,0.101706){\color[rgb]{0,0,0}\makebox(0,0)[lb]{\smash{$u = 0$}}}%
  \end{picture}%
\endgroup